\documentclass[10pt]{article}

\usepackage{amsmath, amsfonts, amsthm, amssymb, amsrefs}
\usepackage{mathtools}
\usepackage{tikz, graphicx, pgfplots, color, xcolor}
	\pgfplotsset{compat=1.12} % indicates version of pgfplots used for consistency
	\usetikzlibrary{positioning}
\usepackage{enumitem}
\usepackage[normalem]{ulem}	% allows sout (strikeout) of text

\usepackage{marginnote}

\usepackage[left=1.2in,right=1.2in,top=1in,foot=1in]{geometry}
\allowdisplaybreaks		% in mathmode environment
\usepackage[linkbordercolor=cyan, pdfborder={0 0 0.5}]{hyperref}
\usepackage{amsthm}
	\newtheorem{thm}{Theorem}[section]
	\newtheorem{lem}[thm]{Lemma}
	
	\newtheorem{cor}[thm]{Corollary}
	\theoremstyle{definition}
		\newtheorem{defn}[thm]{Definition}
		\newtheorem{ex}[thm]{Example}
	\newtheoremstyle{TheoremNum}
        {\topsep}{\topsep}              %%% space between body and thm
        {\itshape}                      %%% Thm body font
        {}                              %%% Indent amount (empty = no indent)
        {\bfseries}                     %%% Thm head font
        {.}                             %%% Punctuation after thm head
        { }                             %%% Space after thm head
        {\thmname{#1}\thmnote{ \bfseries #3}}%%% Thm head spec
    \theoremstyle{TheoremNum}
	\newtheorem{thmref}{Theorem}

	\newcommand{\rmk}{\noindent {\bf Remark.\quad}}

	\newtheorem{prb}{Problem}

	\newcommand{\df}{\bf\emph}		% definitions within definition (bold)

\newcommand{\eq}[1]{\begin{align*}#1\end{align*}}
\newcommand{\eqn}[1]{\begin{align}#1\end{align}}
\newcommand{\ds}{\displaystyle}			% \displaystyle
\renewcommand{\emptyset}{\varnothing}	% varemptyset
\renewcommand{\epsilon}{\varepsilon}	% varepsilon
\renewcommand{\phi}{\varphi}			% varphi
\renewcommand{\tilde}[1]{\widetilde{#1}}% tilde wide
\definecolor{orange}{RGB}{250, 140, 0}
	
\definecolor{turq}{RGB}{0, 160, 160}

\newcommand\mc[1]{\mathcal{#1}}

\newcommand{\rr}{\ensuremath{\mathbb{R}}}

\usetikzlibrary{cd, arrows, automata}
	\newcommand*{\arrow}[1]{\arrow[#1]}

\usepackage{chemfig}
\newcommand{\tikzc}[1]{\begin{center}\begin{tikzpicture}#1\end{tikzpicture}\end{center}}

\DeclareMathOperator{\Ker}{ker}			% kernel
\DeclareMathOperator{\Img}{im}			% image
\DeclareMathOperator{\Span}{span}		% linear span
\DeclareMathOperator{\codim}{codim}		% codim
\DeclareMathOperator{\sgn}{sgn}			% sign function
\newcommand{\braket}[2]{\langle{#1},{#2}\rangle}	% inner product
		\DeclareFontFamily{U}{mathx}{\hyphenchar\font45}
		\DeclareFontShape{U}{mathx}{m}{n}{
	      <5> <6> <7> <8> <9> <10>
	      <10.95> <12> <14.4> <17.28> <20.74> <24.88>
	      mathx10
	      }{}
		\DeclareSymbolFont{mathx}{U}{mathx}{m}{n}
		\DeclareFontSubstitution{U}{mathx}{m}{n}
		\DeclareMathAccent{\widecheck}{0}{mathx}{"71}
\newcommand{\cl}[1]{\overline{#1}}		% top. closure
\newcommand{\union}{\cup}				% union
\newcommand{\bdy}{\partial}				% boundary
\newcommand{\intr}[1]{{#1}^0}			% interior
\newcommand{\kk}{\kappa}

\newcommand{\bb}[1]{{\boldsymbol{#1}}}

\newcommand{\rrp}{\rr_{\geq 0}}
\newcommand{\rrpp}{\rr_{>0}}
\newcommand{\FR}{\ensuremath{\rightarrow}}

\title{A generalization of Birch's theorem and vertex-balanced steady states for generalized mass-action systems}
\date{}

\author{
   Gheorghe Craciun \thanks{Department of Mathematics and Department of Biomolecular Chemistry, University of Wisconsin--Madison}
\and
Stefan M\"uller \thanks{Faculty of Mathematics, University of Vienna}
\and
Casian Pantea \thanks{Department of Mathematics, West Virginia University}
\and
Polly Y. Yu \thanks{Department of Mathematics, University of Wisconsin--Madison}
}

\begin{document}
%% adjust spacing above/below displaystyle/eqn mode.
%	\abovedisplayskip=10pt
%	\belowdisplayskip=10pt

\maketitle
%\mytoc

%%%%%%%%%%%%%%%%%%%%%%%%%%%%%%%%%%%%%%%%%%%%%%%%%%%%%%%%%%%%%%%%%%%%%%%%%%%%%%%%%%%%%%%%%%%%%%%%%%%%%%%%%%%%%%%%%%%%%%%

\begin{abstract}
Mass-action kinetics and its generalizations appear in mathematical models of (bio-)\-chemi\-cal reaction networks, population dynamics, and epidemiology. The dynamical systems arising from directed graphs are generally non-linear and difficult to analyze. One approach to studying them is to find conditions on the network which either imply or preclude certain dynamical properties. For example, a \emph{vertex-balanced steady state} for a \emph{generalized} mass-action system is a state where the net flux through every vertex of the graph is zero. In particular, such steady states admit a monomial parametrization. 
The problem of existence and uniqueness of vertex-balanced steady states can be reformulated in two different ways, one of which is related to \emph{Birch's theorem} in statistics, and the other one to the bijectivity of generalized polynomial maps, similar to maps appearing in geometric modelling. We present a generalization of Birch's theorem, by providing a \emph{sufficient condition} for the existence and uniqueness of vertex-balanced steady states. \\

\end{abstract}

\section{Introduction}
\label{sec:Intro}

	Reaction networks are commonly used to model natural phenomena in disciplines ranging from chemistry, biochemistry, epidemiology to population dynamics. In these systems, entities interact to form other entities as prescribed by a directed graph, the reaction network. For example, the reaction network
	\eq{
	%\label{eq:introEnzyme}
	&\schemestart
		
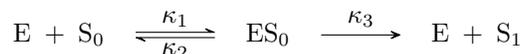
\chemfig{E \+ S_0\, }
			\arrow{<=>[$\kk_1$][$\kk_2$]}
		\chemfig{\, ES_0 \,}
			\arrow{->[$\kk_3$]}
		\chemfig{\, E \+ S_1}
	\schemestop
	}
describes an enzymatic system, where a substrate $\rm{S}_0$ is converted into a product $\rm{S}_1$ by an enzyme $\rm{E}$ via an intermediate species $\rm{ES}_0$. 
	
%	The rate at which these interactions occur is often assumed to be of a specific form; for example, reaction networks in chemistry and biochemistry ...  mass-action kinetics ... The resulting mathematical model is a system of ordinary differential equations. 
	The concentrations of the chemical species in a network are often modelled by a system of ordinary differential equations. One of the most common assumptions in chemistry and biochemistry is that of \emph{mass-action kinetics}, where the reaction rate is proportional to the concentrations of its reactants. According to mass-action kinetics, the reaction $\rm{E + S}_0 \FR \rm{ES}_0$ proceeds at rate $\kk_1 [\rm{E}][\rm{S}_0]$, where $\kk_1 > 0$ is a \emph{rate constant}, and $[\chemfig{X}]$ is the concentration of species $\chemfig{X}$ as a function of time $t$. The rates of change of the concentrations of $\chemfig{E}$, $\chemfig{S_0}$ and $\chemfig{ES_0}$ \emph{due to this single reaction} are
	\eq{
		- \frac{d[\chemfig{E}]}{dt}
		= - \frac{d[\chemfig{S_0}]}{dt}
		= \frac{d[\chemfig{ES_0}]}{dt}
		= \kk_1 [\chemfig{E}][\chemfig{S_0}].
	}
The rates of change \emph{due to all reactions} in the network is the sum over its individual reactions, e.g.
	\eq{
		\frac{d[\chemfig{E}]}{dt}
		&= - \kk_1 [\chemfig{E}][\chemfig{S_0}]
		+ \kk_2 [\chemfig{ES_0}]
		+ \kk_3 [\chemfig{ES_0}],
		\\ %%%
		\frac{d[\chemfig{S_0}]}{dt}
		&= -\kk_1 [\chemfig{E}][\chemfig{S_0}] + \kk_2 [\chemfig{ES_0}],
		\\ %%%
		\frac{d[\chemfig{ES_0}]}{dt}
		&= \hphantom{-}\kk_1 [\chemfig{E}][\chemfig{S_0}]
		- \kk_2 [\chemfig{ES_0}]
		- \kk_3 [\chemfig{ES_0}],
		\\ %%%
		\frac{d[\chemfig{S_1}]}{dt}
		&= \hphantom{- \kk_1 [\chemfig{E}][\chemfig{S_0}]
		- \kk_2 [\chemfig{ES_0}] 11. }
		\kk_3 [\chemfig{ES_0}].
	}

		Mass-action systems have been studied extensively. Reaction network theory, as initially developed by Horn, Jackson and Feinberg~\cite{HJ72, Fein72, Horn72}, tries to conclude dynamical properties from simple characteristics of the underlying network. Moreover, as the reaction rate constant is usually obtained empirically and thus subjected to uncertainty, an ideal theoretical result does not depend on the precise values of the rate constants; indeed this is the case for many classical results in reaction network theory. %Much of the early work on mass-action systems can be found in the lecture notes \cite{GunaNts, FeinLectNts}.

	Mass-action kinetics assumes that the system is dilute (having low concentrations) and homogeneous (well-mixed). In the context of systems biology, that is not the typical environment; the cell is typically crowded and highly structured. Various models have been developed to account for this difference.
	
	Biochemical systems theory \cite{Savageau1969,Voit2013} proposes \emph{power-law kinetics}, where the exponents (or kinetic orders) in the reaction rate functions need not follow the stoichiometric coefficients. In the catalysis example above, we may want the concentration of \chemfig{E} to be modelled by the equation
	\eq{
		\frac{d[\chemfig{E}]}{dt}
		= - \kk_1 [\chemfig{E}]^\alpha[\chemfig{S_0}]^\beta
		+ \kk_2 [\chemfig{ES_0}]^\gamma
		+ \kk_3 [\chemfig{ES_0}]^\delta
	}
for some constants $\alpha$, $\beta$, $\gamma$, $\delta > 0$. This is an example of power-law kinetics. Classical mass-action kinetics and power-law kinetics can be incorporated into the framework of \emph{generalized mass-action kinetics} as formulated in \cite{MR14} (based on \cite{MR12}).

Generalized mass-action systems can also be used to study mass-action systems that do not obviously admit nice dynamical properties. This is done by \emph{network translation}, where a mass-action system is rewritten as a generalized mass-action system, where the underlying network has better properties (e.g., weakly reversible)~\cite{Joh14}. 

\bigskip

	Generalized mass-action systems are essentially dynamical systems %arising from a directed graph $G = (V,E)$. The systems of ordinary differential equations take the very general form 
	of the form 
	\eqn{
	\label{eq:intro1}
		\frac{d\bb x}{dt} &= \sum_{i \in I} \kk_i \, \bb x^{\bb u_i} \bb v_i,
	}
where $\kk_i \in \rrpp$, and $\bb u_i$, $\bb v_i \in \rr^n$. (For $\bb x \in \rrpp^n$ and $\bb u \in \rr^n$, we are using the notation $\bb x^{\bb u} = x_1^{u_1} x_2^{u_2} \cdots x_n^{u_n}$). For example, any polynomial dynamical system is of the form (\ref{eq:intro1}). Moreover, many classes of nonlinear ODEs can be recast as generalized mass-action systems \cite{Brenig1, Brenig2}.
For a complete definition of generalized mass-action systems, see Section~\ref{sec:Gmak}. In this work, we are interested in the existence and uniqueness of steady states of these systems, as it relates to   geometric properties of the vectors $\{\bb u_i, \bb v_i\}_{i \in I}$.

	In classical mass-action systems, some classes of positive steady states enjoy certain algebraic and dynamical properties. Dating back to Boltzmann's kinetic theory, \emph{detailed-balanced equilibria} can be regarded as thermodynamic equilibria. Their generalization, \emph{complex-balanced equilibria}, are dynamically stable because of the existence of an associated  Lyapunov function~\cite{HJ72, FeinLectNts,GunaNts}, and admit monomial parametrizations~\cite{Toric09}. 
	
	For generalized mass-action systems, the analogue of complex-balanced equilibria are the \emph{vertex-balanced steady states}. Unsurprisingly, the theory of vertex-balanced steady states is quite complicated. Some necessary conditions for stability have been found recently~\cite{BorosMuellerRegensburger2019}. Also, they admit monomial parametrizations that may be very useful in applications~\cite{MR14}.

	In this paper, we are interested in how many (if any) vertex-balanced steady states there are within each invariant affine subspace of a generalized mass-action system. In particular, we aim to understand which reaction networks admit vertex-balanced steady states, and whether they are unique. Interestingly, this question can be reformulated in two different ways, one related to a generalization of Birch's theorem in statistics~\cite{SturmfelsPachter2005}, and the other to the bijectivity of generalized polynomial maps, similar to ones which appear in geometric modelling~\cite{MR12, geommodel}. Indeed, the following questions are essentially equivalent:
	\begin{enumerate}
	\item
		When does a generalized mass-action system have exactly one vertex-balanced steady state within each invariant affine subspace, for any choice of rate constants?
	\item
		Given vector subspaces $S$, $\tilde S \subseteq \rr^n$, when does the intersection\footnote{Thereby, $\bb x \circ \bb y$ denotes the component-wise product of the vectors $\bb x$ and $\bb y$; see Section~\ref{sec:Notation}.}
        $(\bb x_0 + S)\cap (\bb x^* \circ \exp\tilde S^\perp)$ consist of exactly one point, for any $\bb x_0$, $\bb x^*\in \rrpp^n$? 	
    \item
		Given vectors $\bb w^1, \dots, \bb w^n, \tilde{\bb w}^1, \dots, \tilde{\bb w}^n \in \rr^d$, when is the generalized polynomial map on $\rrpp^d$ defined by
		\eq{
			f_{\bb x^*}(\xi) = \sum_{i=1}^n x_i^* \xi^{\tilde{\bb w}^i} \bb w^i
		}
		 bijective onto the relative interior of the polyhedral cone generated by $\bb w^1,\dots, \bb w^n$, for any $\bb x^*\in \rrpp^n$?
	\end{enumerate}
These questions will be expanded upon and explained in detail in Section~\ref{sec:equivform}.

%%%%%%%%%%%%%%%%%%%%%%%%%%%%%%%%%%%

Among the questions above, we initially focus on question 2, which is strongly related to Birch's theorem. One way to state Birch's theorem is: given a vector subspace $S \subseteq \rr^n$, the intersection $(\bb x_0 + S) \cap (\bb x^* \circ \exp S^\perp)$ consists of exactly one point, for any $\bb x_0$, $\bb x^* \in \rrpp^n$. In question 2, we have two vector subspaces $S$, $\tilde S$, so it should not come as a surprise that an additional hypothesis is needed, in order for this intersection to consist of exactly one point.

 This additional hypothesis is given in terms of sign vectors. For a subset $S\in \rr^n$, its set of \emph{sign vectors} $\sigma(S)$ is the image of vectors in $S$ under the coordinate-wise sign function. Its \emph{closure} $\cl{\sigma(S)}$ contains $\sigma(S)$ and all sign vectors where one or more coordinates may be replaced by zeros (see Definition~\ref{def:signvectors}). 
 
 \medskip
 
	One of our main results is the following generalization of Birch's theorem:
\begin{thmref}[\ref{thm:main}]
	Let $S$, $\tilde S \subseteq \rr^n$ be vector subspaces of equal dimension with $\sigma(S) \subseteq \cl{\sigma(\tilde S)}$. Then for any positive vectors $\bb x_0$, $\bb x^* \in \rrpp^n$, the intersection $(\bb x_0 + S) \cap (\bb x^* \circ \exp \tilde S^\perp)$ consists of exactly one point.
\end{thmref}
\noindent
By using this theorem, we obtain a sufficient condition for question 1  in Theorem~\ref{thm:mainCRN}. More precisely, provided that certain conditions hold, we show that if a generalized mass-action system has at least one vertex-balanced steady state, then there is exactly one vertex-balanced steady state within every invariant affine subspace.

	We introduce generalized mass-action systems and vertex-balanced steady states in Section~\ref{sec:Gmak} and \ref{sec:VertexBalEqm} respectively. We prove Theorem~\ref{thm:main} and Theorem~\ref{thm:mainCRN}  in Section~\ref{sec:Result}, and conclude with an example in Section~\ref{sec:MainEx}.

	\subsection{Notation}
	\label{sec:Notation}

	There are several component-wise operations on vectors and matrices that will appear frequently. In the list below, let $\bb x$, $\bb z \in \rr^n$ with $\bb x = (x_1,x_2,\dots, x_n)^T$ and $\bb z = (z_1, z_2, \dots, z_n)^T$. Let $Y = (\bb y_1, \bb y_2, \cdots, \bb y_m)$ be a $n\times m$ matrix. 
	
		We write $\bb x \geq 0$ to mean that every component of the vector is non-negative. Similarly, $\bb x > 0$ means that every component of the vector is positive. We let $\rrp^n = \{ \bb x \in \rr^n : \bb x \geq 0\}$, and $\rrpp^n = \{ \bb x \in \rr^n : \bb x > 0\}$. We denote the cardinality of a set $M$ as $|M|$.
		
		The vector and matrix operations we will use are:
		\eq{
			\bb x^{\bb z} &= \prod_{i=1}^n x_i^{z_i}, \text{ where } \bb x > 0;
			\\
			\bb x^Y &= %\begin{pmatrix} \bb x^{\bb y_1} & \bb x^{\bb y_2} & \ldots & \bb x^{\bb y_m} \end{pmatrix}^T;
			(\bb x^{\bb y_1}, \bb x^{\bb y_2}, \ldots, \bb x^{\bb y_m} )^T, \text{ where } \bb x > 0;
			\\
			\bb x\circ \bb z &= %\begin{pmatrix} x_1z_1 & x_2z_2 & \ldots & x_nz_n \end{pmatrix}^T;
			(x_1z_1, x_2z_2,\ldots, x_nz_n)^T;
			\\
			\frac{\bb x}{\bb z} &= %\begin{pmatrix} \frac{x_1}{z_1} & \frac{x_2}{z_2} & \ldots & \frac{x_n}{z_n} \end{pmatrix}^T;
			\left(\frac{x_1}{z_1}, \frac{x_2}{z_2},\ldots, \frac{x_n}{z_n} \right)^T, \text{ where } \bb z > 0;
			\\
			\exp \bb x &= %\begin{pmatrix} e^{x_1} & e^{x_2} & \ldots & e^{x_n} \end{pmatrix}^T;
			(e^{x_1}, e^{x_2},\ldots, e^{x_n})^T;
			\\
			\log \bb x &= %\begin{pmatrix} \log x_1 & \log x_2 & \ldots & \log x_n \end{pmatrix}^T;
			(\log x_1, \log x_2, \ldots, \log x_n )^T, \text{ where } \bb x > 0 .
		}
When the above operations are applied to a subset of $\rr^n$, they are applied to elements of the set. For example, given a set $S \subseteq \rr^n$, we have $\exp(S) = \{ \exp(\bb x) : \bb x \in S\}$, and $\bb x \circ S = \{ \bb x \circ \bb z : \bb z \in S\}$.

%%%%%%%%%%%%%%%%%%%%%%%%%%%%%%%%%%%%%%%%%%%%%%%%%%%%%%%%%%%%%%%%%%%%%%%%%%%%%%%%%%%%%%%%%%%%%%%%%%%%%%%%%%%%%%%%%%%%%%%
\section{Generalized mass-action systems}
\label{sec:Gmak}

	Consider a simple directed graph $G = (V,E)$ and the corresponding weighted digraph $G_{\boldsymbol\kk} = (V,E, \boldsymbol\kk)$ with $\boldsymbol\kk \in \rrpp^{E}$ providing a positive weight for each edge in $E$. Let $V = \{v_1,v_2,\dots, v_m\}$ be the set of vertices. Given an edge $e = v_i \FR v_j \in E$, we call $v_i$ the {\df{source}} of $e$, and $v_j$ its {\df{target}}. Let us denote by $V_s \subseteq V$ the set of \emph{source vertices}, that is, the set of vertices that are sources of some edges. The weight $\kk_e >0$ on the edge $e = v_i \FR v_j$ is called a {\df{rate constant}}, and we refer to the vector $\boldsymbol\kk \in \rrpp^{E}$ as the {\df{vector of rate constants}}, or more simply as the {\df{rate constants}}. Often, we use the indices of the source and target vertices as edge label, i.e., $\kk_{v_i \FR v_j} = \kk_{ij}$.

	Let $\Phi: V \to \rr^n$ be a map assigning to each vertex $v \in V$ a {\df{stoichiometric complex}} $\Phi(v) \in \rr^n$, and let $\tilde \Phi: V_s \to \rr^n$ be another map that assigns to each source vertex $v \in V_s$ a {\df{kinetic-order complex}} $\tilde \Phi(v) \in \rr^n$. An edge $v_i \FR v_j$ is called a {\df{reaction}}, and the vector $\Phi(v_j) - \Phi(v_i)$ is the {\df{reaction vector}} associated to the edge $v_i \FR v_j$. For convenience, we often write $\bb y_i$ instead of $\Phi(v_i)$, and $\tilde{\bb y}_i$ instead of $\tilde \Phi(v_i)$. The graph $G$ and the two maps $\Phi$, $\tilde \Phi$ on $G$ provide all the ingredients needed to define a generalized reaction network, while the weighted digraph $G_{\boldsymbol\kk}$ and the maps $\Phi$, $\tilde \Phi$ are all that is needed to define a generalized mass-action system.

	\begin{defn}
		A {\df{generalized reaction network}} is given by $(G, \Phi, \tilde \Phi)$, where $G = (V,E)$ is a simple directed graph, and $\Phi: V \to \rr^n$, $\tilde \Phi: V_s \to \rr^n$ respectively assign to each vertex a {\df{stoichiometric complex}} and to each source vertex a {\df{kinetic-order complex}}.
	\end{defn}

	\rmk We follow the definition of a generalized reaction network given by M\"uller and Regensburger in \cite{MR14}, rather than the one given in \cite{MR12}. In particular, we do not assume that the maps $\Phi$ and $\tilde \Phi$ are injective.\\
	
	\noindent
	\rmk Throughout this paper, we are concerned with generalized reaction networks where $V_s = V$. The digraphs $\Phi(G)$ and $\tilde \Phi(G)$ are two {\emph{Euclidean embedded graphs}}~\cite{PolylToric, Brunner, RevCY}.
	One of the equivalent definitions of a (classical) \emph{reaction network} is a directed graph $\mc G = (\mc V, \mc E)$, where the set $\mc V$ of vertices (complexes) is a subset of $\rr^n$. 
	Using the notation above, a reaction network is given by $\mc G = \Phi(G)$, where $\Phi$ is injective~\cite{MR14}.

\begin{figure}[h]
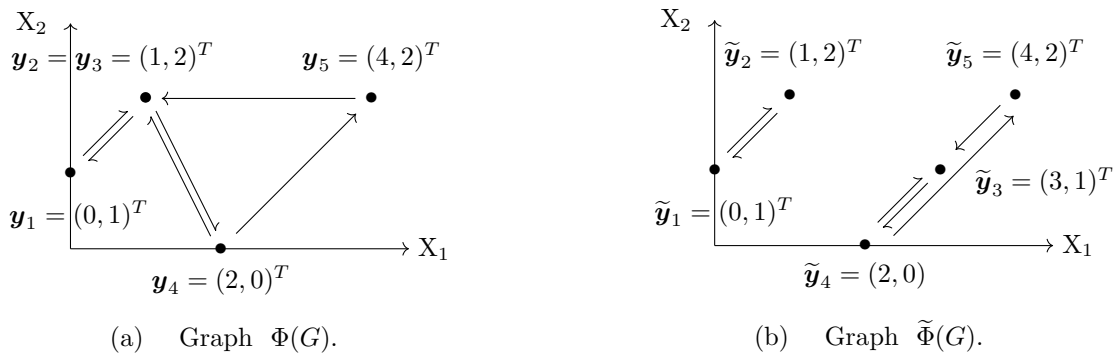

	\centering
		\tikz{
		\node at (+0.4,-1.2) [right]{(a) \quad {Graph } $\Phi(G)$.};
			\draw [->] (0,0)--(4.5,0) node [right] {$\rm{X}_1$};
			\draw [->] (0,0)--(0,3) node [left = 5 pt] {$\rm{X}_2$};
			\node (1) at (0,1) {$\bullet$};
			\node (2) at (1,2) {$\bullet$};
			\node (3) at (1,2) {$\bullet$};
			%\draw  (1,2.02) circle (4pt);
			\node (4) at (2,0) {$\bullet$};
			\node (5) at (4,2) {$\bullet$};
			\node [below =1pt of 1] {\hspace{0.22cm}$\bb y_1 = (0,1)^T$};
			\node [above =1pt of 2] {$\bb y_2 = \bb y_3 = (1,2)^T$\hspace{0.82cm}\,};
			\node [below =-3pt of 4] {$\bb y_4 = (2,0)^T$};
			\node [above = 1pt of 5] {$\bb y_5 = (4,2)^T$};
		\draw [{->[harpoon]}, transform canvas={xshift=-1pt, yshift=1pt}] (1) -- (2) ;
		\draw [{->[harpoon]}, transform canvas={xshift=1pt, yshift=-1pt}] (2) -- (1) ;
		\draw [{->[harpoon]}, transform canvas={xshift=1pt, yshift=1pt}] (3) -- (4) ;
		\draw [{->[harpoon]}, transform canvas={xshift=-1pt, yshift=-1pt}] (4) -- (3) ;
		\draw [->] (4) -- (5);
		\draw [->] (5) -- (3);
		}
		\hspace{2cm}
		\tikz{
		\node at (+0.4,-1.2) [right]{(b) \quad {Graph } $\tilde \Phi(G)$.};
			\draw [->] (0,0)--(4.5,0) node [right] {$\rm{X}_1$};
			\draw [->] (0,0)--(0,3) node [left  = 5 pt] {$\rm{X}_2$};
			\node (1) at (0,1) {$\bullet$};
			\node (2) at (1,2) {$\bullet$};
			\node (3) at (3,1) {$\bullet$};
			\node (4) at (2,0) {$\bullet$};
			\node (5) at (4,2) {$\bullet$};
			\node [below =1pt of 1] {\hspace{0.25cm}$\tilde{\bb y}_1 = (0,1)^T$};
			\node [above =1pt of 2] {\hspace{0.1cm}$\tilde{\bb y}_2 =(1,2)^T$};
			\node [below right =-15pt of 3] {\hspace{0.5cm}$\tilde{\bb y}_3 =(3,1)^T$};
			\node [below =-3pt of 4] {$\tilde{\bb y}_4 = (2,0)$};
			\node [above = 1pt of 5] {$\tilde{\bb y}_5 = (4,2)^T$};
		\draw [{->[harpoon]}, transform canvas={xshift=-1pt, yshift=1pt}] (1) -- (2) ;
		\draw [{->[harpoon]}, transform canvas={xshift=1pt, yshift=-1pt}] (2) -- (1) ;
		\draw [{->[harpoon]}, transform canvas={xshift=1pt, yshift=1pt}] (3) -- (4) ;
		\draw [{->[harpoon]}, transform canvas={xshift=-1pt, yshift=3pt}] (4) -- (3) ;
		\draw [->, transform canvas={xshift=5pt, yshift=-2pt}] (4)--(5);
		\draw [->, transform canvas={xshift=0pt, yshift=2pt}] (5) -- (3);
		}
		\hspace{0cm}\vspace{0.25cm}
		\caption{
		(a) the stoichiometric complex map $\Phi(v_i) = \bb y_i$, (b) the kinetic-order complex map $\tilde \Phi(v_i) = \tilde{\bb y}_i$ (both from $V$ to $\rr^2$),
		and the resulting Euclidean embedded graphs $\Phi(G)$ and $\tilde \Phi(G)$.
		% The vertex $v_3$ is mapped differently by $\Phi$ and $\tilde \Phi$, in such a way that the number of connected components and number of vertices are different in the graphs $\Phi(G)$ and $\tilde \Phi(G)$.
		}
		\label{fig:ExampleEmbedding}
	\end{figure}

\begin{ex}
\label{ex:intro}
	To illustrate the terminology above, we consider a directed graph $G = (V,E)$ and the corresponding weighted digraph $G_{\boldsymbol\kk} = (V,E,\boldsymbol\kk)$:
\end{ex}
	\tikzc{
		\node (1) at (0,0) {$\bullet$};
		\node (2) at (2,0) {$\bullet$};
		\node (3) at (5,0) {$\bullet$};
		\node (5) at (7,0) {$\bullet$};
		\node (4) at (6, -1.4) {$\bullet$};
		\node [left=1pt of 1] {$v_1$};
		\node [right=1pt of 2] {$v_2$};
		\node [left=1pt of 3] {$v_3$};
		\node [below=1pt of 4] {$v_4$};
		\node [right=1pt of 5] {$v_5$};
		\draw [{->[harpoon]}, transform canvas={yshift=1.5pt}] (1) -- (2) node [midway, above] {$\kk_{12}$};
		\draw [{->[harpoon]}, transform canvas={yshift=-1.5pt}] (2) -- (1) node [midway, below] {$\kk_{21}$};
		\draw [{->[harpoon]}, transform canvas={xshift=1.5pt}] (3) -- (4) node [midway, right] {$\kk_{34}$};
		\draw [{->[harpoon]}, transform canvas={xshift=-1.5pt}] (4) -- (3) node [midway, left] {$\kk_{43}$};
		\draw [->] (4) -- (5) node [midway, right] {$\kk_{45}$};
		\draw [->] (5) -- (3) node [midway, above] {$\kk_{53}$};
	}
The set of vertices is $V = \{ v_1,v_2,v_3, v_4, v_5\}$, which coincides with the set of source vertices $V_s$. The set of edges is $E = \{v_1 \FR v_2, \, v_2 \FR v_1, \, v_3 \FR v_4, \, v_4 \FR v_3,\, v_4 \FR v_5,\, v_5 \FR v_3 \}$. The maps $\Phi$ and $\tilde \Phi$ (both from $V$ to $\rr^2$) are given in Figure~\ref{fig:ExampleEmbedding}.

Note that the vertex $v_3$ is mapped differently by $\Phi$ and $\tilde \Phi$. Indeed, $v_3$ is mapped by $\Phi$ to the stoichiometric complex $(1,2)^T$ and by $\tilde \Phi$ to the kinetic-order complex $(3,1)^T$.
Further, note that $\Phi$ maps the vertices $v_2$ and $v_3$ to the same stoichiometric complex, whereas $\tilde \Phi$ maps $v_2$ and $v_3$ to different kinetic-order complexes. Hence, the number of connected components and the number of vertices are different in the graphs $\Phi(G)$ and $\tilde \Phi(G)$.

%\noindent
%\begin{minipage}{\textwidth}
%	There are many reaction networks that can give rise to the generalized reaction network above. For example, we may assume power-law kinetics on the following reaction network:\\
%
%	\tikzc{
%		\schemestart
%			\chemfig{X_2}
%				\arrow(1--2){<=>}
%			\chemfig{X_1\+2X_2 }
%				\arrow(@2--4){<-}
%			\chemfig{4X_1\+2X_2}
%				\arrow(@4--5){<-}[230,]
%			\chemfig{2X_1}
%				\arrow(@5--@2){<=>}
%		\schemestop
%	}
%\end{minipage}

\bigskip

	Now we are in a position to define generalized mass-action systems and the associated dynamical systems.

	\begin{defn}
		A {\df{generalized mass-action system}} is given by $(G_{\boldsymbol\kk}, \Phi, \tilde \Phi)$, where $(G,\Phi,\tilde \Phi)$ is a generalized reaction network, with directed graph $G = (V,E)$, and $\boldsymbol\kk \in \rrpp^{E}$ is a vector of rate constants.
	\end{defn}
	\begin{defn}
	For a generalized mass-action system $(G_{\boldsymbol\kk}, \Phi, \tilde \Phi)$, the {\df{associated dynamical system}} on $\rrpp^n$ is given by
		\eqn{
		\label{eq:gmasode}
			\frac{d\bb x}{dt}
			&= \sum_{v_i \FR v_j \in E} \kk_{ij} \bb x^{\tilde{\bb y}_i} (\bb y_j - \bb y_i ).
		}
	\end{defn}

	As the ODE system (\ref{eq:gmasode}) is our main object of interest, we pause to make two observations. First, the rate of change $\frac{d\bb x}{dt}$ is restricted to the {\df{stoichiometric subspace}} $S = \Span_\rr \{ \bb y_j - \bb y_i : v_i \FR v_j \in E\}$. Consequently, every trajectory $\bb x(t)$ of this dynamical system is restricted to a {\df{stoichiometric compatibility class}} $\bb x(0) + S$. Second, if $v_i \FR v_j$ is a reaction and $\Phi(v_i) = \Phi(v_j)$, then this particular reaction does not contribute to the dynamics.%
%	\footnote{
%		In the definition of a generalized reaction network in~\cite{MR14}, the assumption that the graph $G$ has no self-loops is superfluous, because a self-loop contributes a zero reaction vector and does not affect the dynamics.
%	}%%%%

	It is sometimes more convenient to write the ODE system (\ref{eq:gmasode}) in matrix form. Let $Y \in \rr^{n\times m}$ be the {\df{stoichiometric complex matrix}}, the $j$-th column of which is the stoichiometric complex $\bb y_j$. Let the {\df{kinetic-order complex matrix}} $\tilde Y \in \rr^{n\times m}$ be defined analogously; in particular, its $j$-th column is the kinetic-order complex $\bb{\tilde y}_j$ if $v_j \in V_s$ and $\bb 0$ if $v_j \not\in V_s$.\footnote{%%%
The choice of $\tilde{\bb y}_j = \bb 0$ when $v_j\not\in V_s$ is arbitrary, since the $j$-th column of $\tilde Y$ does not appear in the equations that are of interest to us~\cite{MR14}. In particular, it does not affect the vector $A_{\boldsymbol\kk} \bb x^{\tilde Y}$ and hence does not contribute to the right-hand side $YA_{\boldsymbol\kk} \bb x^{\tilde Y}$ of the system of differential equations (\ref{eq:GmasMtx}).
}
	%%%%%
	 Let $A_{\boldsymbol\kk} \in \rr^{m\times m}$ be the negative transpose of the Laplacian of the weighted directed graph $G_{\boldsymbol\kk}$, i.e.,
	\eqn{
	\label{eq:Laplacian}
		(A_{\boldsymbol\kk})_{ij} &= \left\{
			\begin{array}{cc}\ds
				\vphantom{\sum_{y_i}} \kk_{ji} 	& \text{if } v_j \FR v_i \in E , \\ \ds
				 - \sum_{v_j \FR v_k \in E} \kk_{jk} & \text{if } i = j ,
				\\ \ds
				\vphantom{\sum^{y_i}} 0 & \text{otherwise}.
			\end{array}
		\right.
	}
The dynamical system (\ref{eq:gmasode}) can be rewritten as
	\eqn{
	\label{eq:GmasMtx}
		\frac{d\bb x}{dt} &= Y A_{\boldsymbol\kk} \bb x^{\tilde Y}.
	}

\begin{ex}
	Returning to Example~\ref{ex:intro}, the dynamical system associated to $(G_{\boldsymbol\kk}, \Phi, \tilde \Phi)$ is
	\eq{
	%\label{eq:ExampleODE}
	\begin{aligned}
	%%%
		\frac{d}{dt} \begin{pmatrix} x_1 \\ x_2 \end{pmatrix}
		&= \kk_{12} x_2 \begin{pmatrix} 1 \\ 1 \end{pmatrix}
		+ \kk_{21} x_1x_2^2 \begin{pmatrix} -1 \\ -1 \end{pmatrix}
		+ \kk_{34} x_1^3x_2 \begin{pmatrix} 1 \\ -2 \end{pmatrix}
		%%%%%%%%
        \\
        %%%%%%%%
		& \quad + \, \kk_{43} x_1^2 \begin{pmatrix}-1 \\ 2 \end{pmatrix}
		+ \kk_{45} x_1^2 \begin{pmatrix} 2 \\ 2 \end{pmatrix}
		+ \kk_{53} x_1^4 x_2^2 \begin{pmatrix} -3 \\ 0 \end{pmatrix},
		%%%
		\end{aligned}
	}
where each term corresponds to an edge in the graph $G$. Expanding the equations, we recognize it to be a polynomial (more generally, a power-law) dynamical system:
	\eqn{
	\label{eq:exintro}
		\begin{aligned}
		\ds \frac{dx_1}{dt} &= \kk_{12} x_2 - \kk_{21} x_1x_2^2 + \hphantom{2} \kk_{34} x_1^3x_2 -\hphantom{2} \kk_{43} x_1^2 + 2 \kk_{45} x_1^2 - 3 \kk_{53} x_1^4x_2^2,
		\\
		\ds \frac{dx_2}{dt} &=\kk_{12} x_2 - \kk_{21} x_1x_2^2 - 2\kk_{34} x_1^3x_2 + 2\kk_{43} x_1^2 + 2 \kk_{45} x_1^2.
		\end{aligned}
	}
Its stoichiometric subspace is $S = \rr^2$. The stoichiometric complex matrix and kinetic-order complex matrix are
	\eq{
		Y =
			(\bb y_1 , \bb y_2 , \bb y_3 , \bb y_4 , \bb y_5 )
		= \begin{pmatrix}
			0 & 1 & 1 & 2 & 4\\
			1 & 2 & 2 & 0 & 2
		\end{pmatrix}
		%%%
		\quad \text{and} \quad
		%%%
		\tilde Y =
			(\tilde{\bb y}_1 , \tilde{\bb y}_2 , \tilde{\bb y}_3 , \tilde{\bb y}_4 ,  \tilde{\bb y}_5)
		= \begin{pmatrix}
			0 & 1 & 3 & 2 & 4 \\
			1 & 2 & 1 & 0 & 2
		\end{pmatrix} ,
	}
respectively. The matrix
	\eq{
		A_{\boldsymbol\kk} =
		\begin{pmatrix}
			-\kk_{12} & \kk_{21} \\
			\kk_{12} & -\kk_{21} \\
			& & -\kk_{34} & \kk_{43} 		& \kk_{53}\\
			&& \kk_{34} & -\kk_{43} -\kk_{45} &0 \\
			&& 	0	& \kk_{45} & -\kk_{53}
		\end{pmatrix}
	}
is the negative transpose of the Laplacian of the weighted digraph $G_{\boldsymbol\kk}$. %We refer to $(G, \Phi, \tilde \Phi)$ as the \emph{generalized reaction network}, while $(G_\kk, \Phi, \tilde \Phi)$ is the \emph{generalized mass-action system}, and Equation (\ref{eq:ExampleODE}) is the \emph{associated dynamical system} of $(G_\kk, \Phi, \tilde \Phi)$.
\end{ex}
%\setcounter{section}{3}
%\setcounter{thm}{4}
%%%%%%%%%%%%%%%%%

	The definitions above and Example~\ref{ex:intro} are relatively abstract; one may wonder how generalized mass-action systems show up in applications. Suppose we are interested in modelling the following chemical system:
	
	\tikzc{
		\schemestart
			\chemfig{X_2}
				\arrow(1--2){<=>}
			\chemfig{X_1\+2X_2 }
				\arrow(@2--4){<-}
			\chemfig{4X_1\+2X_2}
				\arrow(@4--5){<-}[230,]
			\chemfig{2X_1}
				\arrow(@5--@2){<=>}
		\schemestop
	}
Let us assume that based on experimental data, the reaction rate functions are as shown below, with rate constants $\kk_{ij} >0$:
	\begin{center}
	\begin{tabular}{l|l}
		Reaction & \, \, Rate Function \, \,  \\
		\hline 
		\chemfig{X_2} \FR \, \chemfig{X_1 + 2X_2} & $\kk_{12} \, x_2$ \\
		\chemfig{X_1 + 2X_2} \FR\, \chemfig{X_2} & $\kk_{21} \, x_1x_2^2$ \\
		\chemfig{X_1+2X_2} \FR\, \chemfig{2X_1} & $\kk_{34} \, x_1^3x_2$ \\
		\chemfig{2X_1} \FR\, \chemfig{X_1 + 2X_2} & $\kk_{43} \, x_1^2 $ \\
		\chemfig{2X_1} \FR\, \chemfig{4X_1 + 2X_2} & $\kk_{45} \, x_1^2$ \\
		\chemfig{4X_1+2X_2} \FR\, \chemfig{X_1+2X_2} & $\kk_{53} \, x_1^4x_2^2$
	\end{tabular}
	\end{center}
Note that the third reaction follows power-law kinetics, but \emph{not} classical mass-action kinetics.  Moreover, note that the second and third reactions are an example of branching reactions, that is, two reactions with the same source complex, different target complexes, and (most importantly) with different kinetics: mass-action for the second reaction and power-law for the third.
Exactly this chemical system is specified as a generalized mass-action system $(G_{\boldsymbol \kk}, \Phi,\tilde\Phi)$ in Example~\ref{ex:intro},
see Figure 1(a) for the reaction network and Figure 1(b) for the reaction rate functions.
The system of ordinary differential equations modelling this system is precisely (\ref{eq:exintro}). \\

	\rmk
	We defined a generalized reaction network as a triple $(G, \Phi, \tilde \Phi)$. As pointed out in an earlier remark, if $\Phi$ is injective, then $\Phi(G)$ is a (classical) reaction network. A classical mass-action system can be obtained as a special case of a generalized mass-action system $(G_{\boldsymbol \kk}, \Phi, \tilde \Phi)$, where $\Phi$ is injective and $\tilde \Phi = \left. \Phi \right|_{V_s}$~\cite{MR12}. 
It is thus natural to extend some of the standard definitions for classical mass-action systems to generalized mass-action systems.\\

	We say the underlying graph $G$ is {\df{weakly reversible}} if every connected component of $G$ is strongly connected, i.e., every edge is part of a directed cycle. We have already defined the stoichiometric subspace $S$ as the span of reaction vectors. Whenever $V_s = V$ (in particular, when $G$ is weakly reversible), we define its kinetic analogue, the {\df{kinetic-order subspace}} $\tilde S = \Span_\rr \{ \tilde{\bb y}_j - \tilde{\bb y}_i : v_i \FR v_j \in E\}$.

	The {\df{stoichiometric deficiency}} of the generalized reaction network $(G, \Phi, \tilde \Phi)$ is the non-negative integer
	\eqn{
		\delta_G = |V| - \ell_G - \dim S,
	}
where $|V|$ is the number of vertices in $G$, $\ell_G$ is the number of connected components of $G$, and $S$ is the stoichiometric subspace. From the equivalent definition $\delta_G = \dim(\Ker Y \cap \Img I_E)$, where $I_E$ is the incidence matrix of $G$, it follows that $\delta_G$ is a non-negative integer~\cite{Joh14}. In the case when $G$ is weakly reversible, we also have the formula~\cite{GunaNts}
	\eqn{
		\delta_G = \dim (\Ker Y \cap \Img A_{\boldsymbol\kk}).
	}
Whenever $V_s = V$, the {\df{kinetic-order deficiency}} is defined as the non-negative integer
	\eqn{
		\tilde \delta_G = |V| - \ell_G - \dim \tilde S,
	}
where $\tilde S$ is the kinetic-order subspace. \\

\rmk In the definitions above, $|V|$ is the number of vertices in the underlying abstract graph $G$, not necessarily the number of distinct stoichiometric complexes or the number of kinetic-order complexes; $\ell_G$ is the number of connected components of $G$, not necessarily the number of connected components in $\Phi(G)$ or $\tilde \Phi(G)$.  \\

In Example~\ref{ex:intro}, we have a weakly reversible network with $|V| = 5$ vertices and $\ell_G = 2$ connected components. We already observed that the stoichiometric subspace $S$ is all of $\rr^2$. However, the kinetic-order subspace is $\tilde S = \Span_\rr (1,1)^T$. The stoichiometric deficiency in this example is $\delta_G = 5 - 2 -2 = 1$, but the kinetic-order deficiency is $\tilde \delta_G = 5 - 2 - 1 = 2$. \\

\begin{ex}
\label{ex:futilecycle}
	We have seen earlier that generalized mass-action systems arise naturally from power-law kinetics. This example illustrates how generalized mass-action systems also arise naturally in the study of mass-action systems, via a process called {\df{network translation}}~\cite{Joh14, Johnston2015}. Network translation produces a generalized mass-action system that has the same dynamics as the original mass-action system. We look at the $n$-site distributive phosphorylation-dephosphorylation system under mass-action kinetics.

	This example first appeared in \cite{Joh14}; below we consider a different translation of the same mass-action system. Under the original definition of generalized mass-action system in \cite{MR12}, which requires the stoichiometric complex map $\Phi$ and the kinetic-order complex map $ \tilde \Phi$ to be injective, the translated network presented below would not have been a well-defined generalized reaction network. However, the later definition in \cite{MR14} removes the requirements that $\Phi$ and $\tilde \Phi$ are injective. As a result, many more dynamical systems can be written as a generalized mass-action system, and for this example, a more natural translation exists for the $n$-site distributive phosphorylation-dephosphorylation system.

	Let $\chemfig{E}$, $\chemfig{F}$ be enzymes that catalyze the phosphorylation and dephosphorlyation processes, by forming intermediates $\chemfig{ES_j}$ and $\chemfig{FS_j}$ respectively. The $n$-site distributive phosphorylation-dephosphorylation system consists of the following reactions:
\end{ex}
\noindent
\begin{minipage}{\textwidth}
	\eq{
	&\schemestart
		\chemfig{E \+ S_0 \,}
			\arrow{<=>}
		\chemfig{\, ES_0 \,}
			\arrow{->}
		\chemfig{\, E \+ S_1 \,}
			\arrow{<=>}
		\chemfig{\, ES_1 \,}
			\arrow{->}
		\chemfig{\cdots}
	\schemestop
	\\ %%
	\hphantom{hi}&\hspace{5cm}
	\schemestart
		\chemfig{\cdots}
			\arrow{->}
		\chemfig{\, E\+S_{n-1} \, }
			\arrow{<=>}
		\chemfig{\, ES_{n-1}\, }
			\arrow{->}
		\chemfig{\, E \+ S_n}
	\schemestop
	\\
	&\hphantom{hi} \\ %%%%
	&\schemestart
		\chemfig{F \+ S_0\, }
			\arrow{<-}
		\chemfig{\, FS_1 \,}
			\arrow{<=>}
		\chemfig{\, F \+ S_1 \,}
			\arrow{<-}
		\chemfig{\, FS_2 \,}
			\arrow{<=>}
		\chemfig{\cdots}
	\schemestop
	\\ %%
	\hphantom{hi}&\hspace{5cm}
	\schemestart
		\chemfig{\cdots}
			\arrow{<=>}
		\chemfig{\, F \+ S_{n-1} \,}
			\arrow{<-}
		\chemfig{\, FS_{n\hphantom{-1}} \,}
			\arrow{<=>}
		\chemfig{\, F \+ S_n}
	\schemestop \\
	&
	}
\end{minipage}
We assume that the reaction rates follow classical mass-action kinetics. There are $3n+3$ species involved, so the system of differential equations modelling their concentrations is defined on $\rrpp^{3n+3}$. We create a generalized mass-action system with the same differential equations by network translation. The main step involves changing the stoichiometric complexes: adding enzyme $\chemfig{F}$ to the series of reactions for phosphorylation by $\chemfig{E}$; and adding enzyme $\chemfig{E}$ to the series of reactions for dephosphorylation by $\chemfig{F}$. This process produces a weakly reversible network:
	\tikzc{
%		\node (1l) at (0,0) [left] {$\chemfig{E + F} + \chemfig{S_0}$};
%		\node (1u) at (1,1.5) {$\chemfig{F} + \chemfig{ES_0}$};
%		\node (1d) at (1,-1.5) {$\chemfig{E + FS_1}$};
%		\node (1r) at (2,0) [right]{$\chemfig{E + F} + \chemfig{S_1}$};
		\node (1l) at (0,0) {};
		\node (1u) at (1,1) {};
		\node (1d) at (1,-1) {};
		\node (1r) at (2,0) {};
		\node at (5,0) {$\cdots \cdots$};
		\draw [->] (1u)--(1r);
		\draw [->] (1d)--(1l);
		\draw [{->[harpoon]}, transform canvas={xshift=-1pt, yshift=1pt}] (1l)--(1u);
		\draw [{->[harpoon]}, transform canvas={xshift=1pt, yshift=-1pt}] (1u)--(1l);
		\draw [{->[harpoon]}, transform canvas={xshift=-1pt, yshift=1pt}] (1d)--(1r);
		\draw [{->[harpoon]}, transform canvas={xshift=1pt, yshift=-1pt}] (1r)--(1d);
		%%%%%%%%%%%%
		\node [left =0pt of 1l] {$\chemfig{E + F} + \chemfig{S_0}$};
		\node [right=0pt of 1r] {$\chemfig{E + F} + \chemfig{S_1}$};
		\node [above=0pt of 1u] {$\chemfig{F} + \chemfig{ES_0}$};
		\node [below=0pt of 1d] {$\chemfig{E + FS_1}$};
	}

	To define the generalized mass-action system, we take a more top-down approach, starting from a graph $G$ with $n$ components and $4n$ vertices:
	\tikzc{
		\node (1l) at (0,0) {$\bullet$};
		\node (1u) at (1,1) {$\bullet$};
		\node (1d) at (1,-1) {$\bullet$};
		\node (1r) at (2,0) {$\bullet$};
		\node (2l) at (4,0) {$\bullet$};
		\node (2u) at (5,1) {$\bullet$};
		\node (2d) at (5,-1) {$\bullet$};
		\node (2r) at (6,0) {$\bullet$};
		\node at (8,0) {$\cdots$};
		\node (nl) at (10,0) {$\bullet$};
		\node (nu) at (11,1) {$\bullet$};
		\node (nd) at (11,-1) {$\bullet$};
		\node (nr) at (12,0) {$\bullet$};
		\draw [->] (1u)--(1r);
		\draw [->] (1d)--(1l);
		\draw [{->[harpoon]}, transform canvas={xshift=-1pt, yshift=1pt}] (1l)--(1u);
		\draw [{->[harpoon]}, transform canvas={xshift=1pt, yshift=-1pt}] (1u)--(1l);
		\draw [{->[harpoon]}, transform canvas={xshift=-1pt, yshift=1pt}] (1d)--(1r);
		\draw [{->[harpoon]}, transform canvas={xshift=1pt, yshift=-1pt}] (1r)--(1d);
		\draw [->] (2u)--(2r);
		\draw [->] (2d)--(2l);
		\draw [{->[harpoon]}, transform canvas={xshift=-1pt, yshift=1pt}] (2l)--(2u);
		\draw [{->[harpoon]}, transform canvas={xshift=1pt, yshift=-1pt}] (2u)--(2l);
		\draw [{->[harpoon]}, transform canvas={xshift=-1pt, yshift=1pt}] (2d)--(2r);
		\draw [{->[harpoon]}, transform canvas={xshift=1pt, yshift=-1pt}] (2r)--(2d);
		\draw [->] (nu)--(nr);
		\draw [->] (nd)--(nl);
		\draw [{->[harpoon]}, transform canvas={xshift=-1pt, yshift=1pt}] (nl)--(nu);
		\draw [{->[harpoon]}, transform canvas={xshift=1pt, yshift=-1pt}] (nu)--(nl);
		\draw [{->[harpoon]}, transform canvas={xshift=-1pt, yshift=1pt}] (nd)--(nr);
		\draw [{->[harpoon]}, transform canvas={xshift=1pt, yshift=-1pt}] (nr)--(nd);
		%%%%%%%%%%%%
		\node [left =0pt of 1l] {$v_0$};
		\node [right=0pt of 1r] {$v_1'$};
		\node [left=0pt of 2l] {$v_1$};
		\node [right=0pt of 2r] {$v_2'$};
		\node [left=0pt of nl] {$v_{n-1}$};
		\node [right=0pt of nr] {$v_n'$};
		\node [above=0pt of 1u] {$z_0$};
		\node [above=0pt of 2u] {$z_1$};
		\node [above=0pt of nu] {$z_{n-1}$};
		\node [below=0pt of 1d] {$w_1$};
		\node [below=0pt of 2d] {$w_2$};
		\node [below=0pt of nd] {$w_n$};
	}

	Although $\Phi$ and $\tilde \Phi$ map vertices to vectors in $\rr^{3n+3}$, to make this example more readable, we will specify the images of $\Phi$ and $\tilde \Phi$ in terms of formal linear combination of species.

The stoichiometric complexes are
	\eq{
		\Phi(v_j) = \chemfig{E + F} + \chemfig{ S_j},
		\quad
		\Phi(v_j') = \chemfig{E + F} +\chemfig{ S_j},
		\quad
		\Phi(z_j) = \chemfig{F + ES_j},
		\quad
		\Phi(w_j) = \chemfig{E + FS_j},
	}
and the kinetic-order complexes are
	\eq{
		\tilde\Phi(v_j) = \chemfig{E + S_j},
		\quad
		\tilde\Phi(v_j') = \chemfig{F + S_j},
		\quad
		\tilde\Phi(z_j) = \chemfig{ES_j},
		\quad
		\tilde\Phi(w_j) = \chemfig{FS_j}.
	}
Note that the map $\Phi$ is not injective, as $\Phi(v_j) = \Phi(v_j')$ for $1 \leq j \leq n-1$. The image of the graph $G$ under $\Phi$ is connected:

	\tikzc{
		\node (1l) at (0,0) {$\bullet$};
		\node (1u) at (1,1) {$\bullet$};
		\node (1d) at (1,-1) {$\bullet$};
		\node (1r) at (2,0) {$\bullet$};
		\node (2l) at (2,0) {$\bullet$};
		\node (2u) at (3,1) {$\bullet$};
		\node (2d) at (3,-1) {$\bullet$};
		\node (2r) at (4,0) {$\bullet$};
		\node at (5.5,0) {$\cdots$};
		\node (nl) at (7,0) {$\bullet$};
		\node (nu) at (8,1) {$\bullet$};
		\node (nd) at (8,-1) {$\bullet$};
		\node (nr) at (9,0) {$\bullet$};
		\draw [->] (1u)--(1r);
		\draw [->] (1d)--(1l);
		\draw [{->[harpoon]}, transform canvas={xshift=-1pt, yshift=1pt}] (1l)--(1u);
		\draw [{->[harpoon]}, transform canvas={xshift=1pt, yshift=-1pt}] (1u)--(1l);
		\draw [{->[harpoon]}, transform canvas={xshift=-1pt, yshift=1pt}] (1d)--(1r);
		\draw [{->[harpoon]}, transform canvas={xshift=1pt, yshift=-1pt}] (1r)--(1d);
		\draw [->] (2u)--(2r);
		\draw [->] (2d)--(2l);
		\draw [{->[harpoon]}, transform canvas={xshift=-1pt, yshift=1pt}] (2l)--(2u);
		\draw [{->[harpoon]}, transform canvas={xshift=1pt, yshift=-1pt}] (2u)--(2l);
		\draw [{->[harpoon]}, transform canvas={xshift=-1pt, yshift=1pt}] (2d)--(2r);
		\draw [{->[harpoon]}, transform canvas={xshift=1pt, yshift=-1pt}] (2r)--(2d);
		\draw [->] (nu)--(nr);
		\draw [->] (nd)--(nl);
		\draw [{->[harpoon]}, transform canvas={xshift=-1pt, yshift=1pt}] (nl)--(nu);
		\draw [{->[harpoon]}, transform canvas={xshift=1pt, yshift=-1pt}] (nu)--(nl);
		\draw [{->[harpoon]}, transform canvas={xshift=-1pt, yshift=1pt}] (nd)--(nr);
		\draw [{->[harpoon]}, transform canvas={xshift=1pt, yshift=-1pt}] (nr)--(nd);
		%%%%%
		\node (f2u) at (4.75,0.75) {};
		\node (f2d) at (4.75,-0.75) {};
		\draw [{->[harpoon]}, transform canvas={xshift=1pt, yshift=-1pt}] (f2u)--(2r);
		\draw [transform canvas={xshift=-1pt, yshift=1pt}] (2r)--(f2u);
		\draw [->] (f2d)--(2r);
		%%%%%
		\node (fnu) at (6.25,0.75) {};
		\node (fnd) at (6.25,-0.75) {};
		\draw [{->[harpoon]}, transform canvas={xshift=-1pt, yshift=1pt}] (fnd)--(nl);
		\draw [transform canvas={xshift=1pt, yshift=-1pt}] (nl)--(fnd);
		\draw [->] (fnu)--(nl);
		%%%%%
%		\draw (2,0.02) circle (4pt);
%		\draw (4,0.02) circle (4pt);
%		\draw (7,0.02) circle (4pt);
		}
One can check that $\dim S = \dim \tilde S = 3n$; therefore, the stoichiometric deficiency and kinetic-order deficiency are $\delta_G = \tilde \delta_G = (4n) - (n) - (3n) = 0$.

%%%%%%%%%%%%%%%%%%%%%%%%%%%%%%%%%%%%%%%%%%%%%%%%%%%%%%%%%%%%%%%%%%%%%%%%%%%%%%%%%%%%%%%%%%%%%%%%%%%%%%%%%%%%%%%%%%%%%%%
\section{Vertex-balanced steady states}
\label{sec:VertexBalEqm}

	Given the dynamical system associated to a generalized mass-action system $(G_{\boldsymbol\kk}, \Phi, \tilde \Phi)$, written either as
	$\frac{d\bb x}{dt} = \sum_{v_i\FR v_j\in E} \kk_{ij} \bb{x}^{\tilde{\bb y}_i} (\bb y_j - \bb y_i)$
	or in matrix form
	$\frac{d\bb x}{dt} = YA_{\boldsymbol\kk} \bb x^{\tilde Y}$, it is natural to ask how many steady states there are. We define the set of {\df{positive steady states}} as
	\eqn{
	\label{def:eqm}
		E_{\boldsymbol\kk} = \{ \bb x \in \rrpp^n : YA_{\boldsymbol\kk} \bb x^{\tilde Y} = \bb 0 \}.
	}

	For a classical mass-action system, an important subset of positive steady states is the set of {\emph{complex-balanced equilibria}}~\cite{HJ72}, also known as {\emph{complex-balancing equilibria}} or {\emph{vertex-balanced equilibria}}~\cite{CraciunGAC}. Horn and Jackson introduced the idea of complex balancing at equilibrium to generalize the physical assumption of detailed balancing at thermodynamic equilibrium~\cite{HJ72}.

	We illustrate the intuition behind the definition of such a steady state before introducing its analogue for a generalized mass-action system. Consider the graph $G$ of the reaction network, and associate to each edge $v_i \FR v_j$ a reaction rate function $\kk_{ij} \bb x^{\bb y_i}$. A concentration vector $\bb x^* \in \rrpp^n$ is a \emph{complex-balanced equilibrium} of the classical mass-action system if at every vertex $v_i \in V$ of the graph, the sum of incoming fluxes balances the sum of outgoing fluxes, that is, for all $v_i \in V$, 
	\eqn{
	\label{def:VBclassical}
		\sum_{v_j \FR v_i \in E} \kk_{ji} (\bb x^*)^{\bb y_j}
		= \sum_{v_i \FR v_j \in E} \kk_{ij} (\bb x^*)^{\bb y_i}.
	}
This occurs if and only if $A_{\boldsymbol\kk} (\bb x^*)^{Y} = \bb 0$~\cite{HJ72}. Clearly, a complex-balanced equilibrium is a positive solution to a system of polynomial equations. Surprisingly, it is also a positive solution to a system of \emph{binomial} equations~\cite{Toric09}. \\

	For a generalized mass-action system, one can define a vertex-balanced steady state analogously: it is a positive steady state at which the net flux is zero across every vertex of the graph, where the flux is given by generalized mass-action kinetics.
	\begin{defn}
	\label{def:CBeqm}
	The set of {\df{vertex-balanced steady states}} for a generalized mass-action system $(G_{\boldsymbol\kk}, \Phi, \tilde \Phi)$ is the set
	\eqn{
		Z_{\boldsymbol\kk} = \{ \bb x \in \rrpp^n : A_{\boldsymbol\kk} \bb x^{\tilde Y} = \bb 0\}.
	}
	\end{defn}

\noindent
Note that $\bb x^* \in Z_{\boldsymbol \kk}$ if and only if for all $v_i \in V$, 
	\eqn{
		\label{def:VBgen}
		\sum_{v_j \FR v_i \in E} \kk_{ji} (\bb x^*)^{\tilde{\bb y}_j}
		= \sum_{v_i \FR v_j \in E} \kk_{ij} (\bb x^*)^{\tilde{\bb y}_i}.
	}

\rmk What we call vertex-balanced steady state here, is also called \emph{complex balancing equilibrium}~\cite{MR12, MR14} or \emph{generalized complex-balanced steady state}~\cite{Joh14}.  \\

We call such a steady state \emph{vertex-balanced} instead of \emph{complex-balanced} to avoid a subtle point of confusion. 
In the case when $\Phi$ is not injective, the balancing of in-fluxes and out-fluxes occurs at each vertex $v \in V$ of the underlying abstract graph $G$. This in turn implies the balancing of fluxes at each stoichiometric complex $\Phi(v) \in \Phi(V)$; however, the converse is generally false. For example, consider the generalized mass-action system given by the weighted digraph $G_{\boldsymbol\kk}$,
	\tikzc{
		\node (1) at (0,0) {$\bullet$};
		\node (2) at (2,0) {$\bullet$};
		\node (3) at (0,-1) {$\bullet$};
		\node (4) at (2,-1) {$\bullet$};
		\node [left=1pt of 1] {$v_1$};
		\node [left=1pt of 3] {$v_4$};
		\node [right=1pt of 2] {$v_2$};
		\node [right=1pt of 4] {$v_3$};
		\draw [->] (1) --(2) node [midway, above] {$\kk$};
		\draw [->] (4)--(3) node [midway, below] {$\kk$};
	}
and the maps $\Phi(v_1) = \Phi(v_4) = 0 \in \rr^1$, $\Phi(v_2) = \Phi(v_3) = 1$, $\tilde \Phi(v_1) =  0$, and  $\tilde \Phi(v_3) = 1$. The associated dynamical system is $\frac{dx}{dt} = \kk - \kk x$, 
which also arises from the classical mass-action system $\Phi(G_{\boldsymbol\kk})$:
	\tikzc{
		\node (1) at (0,0) {$\bullet$};
		\node (2) at (2,0) {$\bullet$};
		\node [left=1pt of 1] {$0$};
		\node [right=1pt of 2] {$
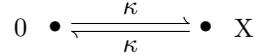
\chemfig{X}$};
		\draw [{->[harpoon]}, transform canvas={yshift=1.5pt}] (1) -- (2) node [midway, above] {$\kk$};
		\draw [{->[harpoon]}, transform canvas={yshift=-1.5pt}] (2) -- (1) node [midway, below] {$\kk$};
	}
% and the maps $\Phi'(v_1)  = 0 \in \rr^1$ and $\Phi'(v_2) = 1$. 
However, the equilibrium $x^* = 1$ is not vertex-balanced for the generalized mass-action system $(G_{\boldsymbol\kk},\Phi,\tilde\Phi)$,
but is complex-balanced for the classical mass-action system $\Phi(G_{\boldsymbol\kk})$.
%\blue{that is, balanced at the level of (stoichiometric) complexes in $G_{\boldsymbol\kk}$.}
\\

\begin{ex}
	To illustrate the definition of vertex-balanced steady states, we consider Example~\ref{ex:intro} again. A vertex-balanced steady state is a point $\bb x = (x_1,x_2)^T \in \rrpp^2$ satisfying five polynomial equations, one equation for each vertex of the graph $G$:
	\eqn{
	\label{eq:exVBss}
		\begin{array}{rrcl}
		v_1:
			& \kk_{12} x_2
			&=&
			\kk_{21} x_1x_2^2, \\
		v_2:
			&\kk_{21} x_1x_2^2
			&=&
			\kk_{12} x_2,	\\
		v_3:
			& \kk_{34} x_1^3x_2
			&=&
			\kk_{43} x_1^2  + \kk_{53} x_1^4x_2^2, \\
		v_4:
			&\quad  (\kk_{43} + \kk_{45}) x_1^2
			&=&
			\kk_{34} x_1^3 x_2, \\
		v_5:
			& \kk_{53} x_1^4x_2^2
			&=&
			\kk_{45} x_1^2. \\
		\end{array}
	}
\end{ex}
However, a positive steady state $\bb x = (x_1,x_2)^T \in \rrpp^2$ has to satisfy only two polynomial equations:
	\eqn{
	\label{eq:exEqm}
	\begin{aligned}
		0 = \frac{dx_1}{dt} &= \kk_{12} x_2 - \kk_{21} x_1x_2^2 + \hphantom{2} \kk_{34} x_1^3x_2 -\hphantom{2} \kk_{43} x_1^2 + 2 \kk_{45} x_1^2 - 3 \kk_{53} x_1^4x_2^2,
		\\
		0 = \frac{dx_2}{dt} &=\kk_{12} x_2 - \kk_{21} x_1x_2^2 - 2\kk_{34} x_1^3x_2 + 2\kk_{43} x_1^2 + 2 \kk_{45} x_1^2.
	\end{aligned}
	}
The two polynomial equations (\ref{eq:exEqm}) are linear combinations of the five polynomial equations (\ref{eq:exVBss}); thus $Z_{\boldsymbol\kk} \subseteq E_{\boldsymbol\kk}$. This follows from the matrix expression of the associated dynamical system $\frac{d\bb x}{dt} = Y( A_{\boldsymbol\kk} \bb x^{\tilde Y})$. \\

	Complex-balanced equilibria of classical mass-action systems have been studied extensively. Some of the classical results extend directly to the case of generalized mass-action systems, even when the maps $\Phi$ and $\tilde \Phi$, assigning stoichiometric complexes and kinetic-order complexes respectively, are not injective. For example, it is known that \cite{MR12, MR14}:~\footnote{%%
Some of these results were first proved in \cite{MR12}, under the assumption that $\Phi$ and $\tilde \Phi$ are injective, but the same proof goes through without these hypotheses.
}%%
	\begin{enumerate}[label=\roman*)]
	\item
		If $Z_{\boldsymbol\kk} \neq \emptyset$ for some ${\boldsymbol\kk} > 0$, then the underlying graph $G$ is weakly reversible.
	\item
		If $Z_{\boldsymbol\kk} \neq \emptyset$ and $\bb x^* \in Z_{\boldsymbol\kk}$, then $Z_{\boldsymbol\kk} = \{ \bb x \in \rrpp^n : \ln \bb x- \ln \bb x^* \in \tilde S^\perp\}
			= \bb x^* \circ \exp \tilde S^\perp$.
	\item
		For a weakly reversible generalized reaction network, $\tilde \delta_G = 0$ if and only if
		$Z_{\boldsymbol \kk} \neq \emptyset$ 
		% the set of vertex-balanced steady states $Z_{\boldsymbol\kk}$ is non-empty 
		for any choices of rate constants ${\boldsymbol\kk} > 0$.
	\item
		For a weakly reversible generalized reaction network, if $\delta_G = 0$, then for any choice of rate constants ${\boldsymbol\kk} > 0$, any positive steady state is a vertex-balanced steady state, i.e., $E_{\boldsymbol\kk} = Z_{\boldsymbol\kk}$.
	\end{enumerate}

	In the example of the $n$-site phosphorlyation-dephosphorlyation system (Example~\ref{ex:futilecycle}), we noted that $\delta_G = \tilde \delta_G = 0$. By statements (iii) and (iv) above, we conclude for any rate constants ${\boldsymbol\kk}$, the set of vertex-balanced steady states $Z_{\boldsymbol\kk}$ is non-empty, and all positive steady states are vertex-balanced. Moreover, the set of positive steady states is given by $E_{\boldsymbol\kk} = Z_{\boldsymbol\kk}  = \bb x^* \circ \exp \tilde{S}^\perp$, where $\bb x^*$ is any positive steady state and $\tilde{S}$ is the kinetic-order subspace, i.e., the vector space spanned by the differences of kinetic-order complexes according to the edges in the graph. It should be noted that the $n$-site phosphorlyation-dephosphorlyation system is \emph{multistationary} when $n \geq 2$~\cite{MAPK1, MAPK2}, i.e., the system admits multiple steady states within the same stoichiometric compatibility class. In other words, for some choices of rate constants, there are multiple vertex-balanced steady states within some stoichiometric compatibility class. This contrasts with a classical complex-balanced mass-action system, where $Z_{\boldsymbol\kk}$ meets every stoichiometric compatibility class at most once.
	
	\begin{ex} 
	   There is another surprising way in which vertex balancing differs from classical complex-balanced mass-action systems. While it is clear that $Z_{\boldsymbol\kk} \subseteq E_{\boldsymbol\kk}$, in generalized mass-action systems it is possible that $\emptyset \neq Z_{\boldsymbol\kk} \subsetneq E_{\boldsymbol\kk}$\footnote{%
	        For classical complex-balanced mass-action systems, it is always the case that $Z_{\boldsymbol\kk} = E_{\boldsymbol\kk}$~\cite{HJ72}.
	    }%
	    . For example, consider the 1-dimensional generalized mass-action system given by the weighted digraph $G_{\boldsymbol\kk}$
	    \tikzc{
		\node (1) at (0,0) {$\bullet$};
		\node (2) at (2,0) {$\bullet$};
		\node (3) at (0,-1.5) {$\bullet$};
		\node (4) at (2,-1.5) {$\bullet$};
		\node [left=1pt of 1] {$v_1$};
		\node [right=1pt of 2] {$v_2$};
		\node [left=1pt of 3] {$v_3$};
		\node [right=1pt of 4] {$v_4$};
		\draw [{->[harpoon]}, transform canvas={yshift=1.5pt}] (1) -- (2) node [midway, above] {$1$};
		\draw [{->[harpoon]}, transform canvas={yshift=-1.5pt}] (2) -- (1) node [midway, below] {$1$};
		\draw [{->[harpoon]}, transform canvas={yshift=1.5pt}] (3) -- (4) node [midway, above] {$5$};
		\draw [{->[harpoon]}, transform canvas={yshift=-1.5pt}] (4) -- (3) node [midway, below] {$5$};
	    }
	    and the maps $\Phi$, $\tilde{\Phi}$ given by
	    \eq{
	        \Phi(v_1) = \tilde{\Phi}(v_1) = 0, 
	        \qquad 
	        \Phi(v_2) = 1, \tilde{\Phi}(v_2) = 3, \\
	        \Phi(v_3) = \tilde{\Phi}(v_3) = 2,
	        \qquad
	        \Phi(v_4) = 3, \tilde{\Phi}(v_4) = 1.
	    }
	    The associated dynamical system is 
	   \eq{
	        \frac{dx}{dt} = 1 - 5x + 5x^2 - x^3.
	   }
	   The point $x^* = 1$ is a vertex-balanced steady state. The system also has two other steady states $x^* \approx 0.27$ and $x^* \approx 3.72$, neither of which satisfy the vertex-balanced condition:
	   \eq{
	    1 = (x^*)^3 \quad \text{and} \quad 
	    5(x^*)^2 = 5 (x^*)^2. 
	   }
	\end{ex}

\bigskip

	In applications, the vector of rate constants ${\boldsymbol\kk} \in \rrpp^E$ is often not known precisely. Surprisingly, some important results for complex-balanced equilibria in classical mass-action systems hold irrespective of the precise values of the rate constants. We are interested in results for vertex-balanced equilibria of generalized mass-action systems that are in this sense independent of the choice of rate constants. We have observed that the solution trajectories are confined to a stoichiometric compatibility class $\bb x_0 + S$, where $\bb x_0 \in \rrpp^n$ is an initial state and $S$ is the stoichiometric subspace. Therefore, our object of study is the intersection $(\bb x_0 + S) \cap Z_{\boldsymbol\kk}$ for any $\bb x_0 \in \rrpp^n$ and any ${\boldsymbol\kk} \in \rrpp^E$.

%%%%%%%%%%%%%%%%%%%%%%%%%%%%%%%%%%%%%%%%%%%%%%%%%%%%%%%%%%%%%%%%%%%%%%%%%%%%%%%%%%%%%%%%%%%%%%%%%%%%%%%%%%%%%%%%%%%%%%%
\section{Problem reformulations}
\label{sec:equivform}

	In the introduction, we have mentioned that the following questions are essentially equivalent:
	\begin{enumerate}
	\item
		When does a generalized mass-action system have exactly one vertex-balanced steady state within each stoichiometric compatibility class, for any choice of rate constants?
	\item
		Given vector subspaces $S$, $\tilde S \subseteq \rr^n$, when does the intersection $(\bb x_0 + S)\cap (\bb x^* \circ \exp\tilde S^\perp)$ consist of exactly one point, for any $\bb x_0$, $\bb x^*\in \rrpp^n$?
	\item
		Given vectors $\bb w^1, \dots, \bb w^n, \tilde{\bb w}^1, \dots, \tilde{\bb w}^n \in \rr^d$, when is the generalized polynomial map on $\rrpp^d$ defined by
		\eq{
			f_{\bb x^*}(\xi) = \sum_{i=1}^n x_i^* \xi^{\tilde{\bb w}^i} \bb w^i
		}
bijective onto the relative interior of the polyhedral cone generated by $\bb w^1,\dots, \bb w^n$, for any $\bb x^*\in \rrpp^n$?
	\end{enumerate}

	Before we discuss the relationship between these problems in detail, let us first make a historical note. When speaking of a weakly reversible classical mass-action system, Horn and Jackson~\cite{HJ72} proved that if the system has at least one complex-balanced equilibrium, then every stoichiometric compatibility class has exactly one complex-balanced equilibrium. Indeed, they showed that every positive steady state of such a system is complex-balanced and locally asymptotically stable within its stoichiometric compatibility class. A complex-balanced equilibrium is \emph{globally} stable within its stoichiometric compatibility class when the network has a single connected component~\cite{DaveGAC}, or is strongly endotactic~\cite{geomGAC}, or when the system is in $\rr^3$~\cite{CasianGAC1, CasianGAC2}. A general proof of global stability of complex-balanced equilibrium within its stoichiometric compatibility class was proposed for all complex-balanced systems in \cite{CraciunGAC}.

\bigskip

	 The first of the three questions above is phrased in the context of reaction networks. We start with a generalized reaction network and suppose that for some rate constants ${\boldsymbol\kk}$, there is a vertex-balanced steady state $\bb x^* \in Z_{\boldsymbol\kk}$. What is a condition (E) on the network $(G, \Phi, \tilde \Phi)$ for the \emph{existence} of a vertex-balanced steady state within every stoichiometric compatibility class? What is a condition (U) on $(G, \Phi, \tilde \Phi)$ so that a vertex-balanced steady state is \emph{unique} within its stoichiometric compatibility class? We would like to obtain conditions for these to hold or fail that are independent of the rate constants ${\boldsymbol\kk}$. More precisely:
	\begin{prb}[1]
	\label{prb:CRNT}
		Let $(G_{\boldsymbol\kk}, \Phi, \tilde \Phi)$ be a generalized mass-action system. Suppose that $Z_{\boldsymbol\kk} \neq \emptyset$. What are conditions (E) and (U) on $(G, \Phi, \tilde \Phi)$, so that the following statements are true?
		\begin{enumerate}\itemsep=-2pt
		\item
		If $(G, \Phi,\tilde\Phi)$ satisfies condition (E), then there is {\bf at least one} vertex-balanced steady state in every stoichiometric compatibility class, i.e., $ (\bb x_0 +S) \cap Z_{\boldsymbol\kk}$ contains at least one point for any $\bb x_0 \in \rrpp^n$.
		\item
		If $(G, \Phi,\tilde\Phi)$ satisfies condition (U), then there is {\bf at most one} vertex-balanced steady state in every stoichiometric compatibility class, i.e., $(\bb x_0 + S) \cap Z_{\boldsymbol\kk}$ contains at most one point for any $\bb x_0 \in \rrpp^n$.
		\end{enumerate}
	\end{prb}

	\bigskip

	Recall that $Z_{\boldsymbol\kk} = \bb x^* \circ \exp \tilde S^\perp$ for any $\bb x^* \in Z_{\boldsymbol\kk}$. Thus, the vertex-balanced steady states within any stoichiometric compatibility class $\bb x_0 +S$ belong to the intersection $(\bb x_0 + S) \cap (\bb x^* \circ \exp \tilde S^\perp)$. This leads us to the following reformulation of Problem~\ref{prb:CRNT}:
	\begin{prb}[2]
	\label{prb:mfld}
		Let $S$, $\tilde S \subseteq \rr^n$ be vector subspaces. What are conditions (E) and (U) on $S$, $\tilde S$, so that the following statements are true?
		\begin{enumerate}\itemsep=-2pt
		\item
		If $S$, $\tilde S$ satisfy condition (E), then $(\bb x_0 + S) \cap (\bb x^* \circ \exp \tilde S^\perp)$ contains {\bf at least one} point, for any $\bb x_0$, $\bb x^* \in \rrpp^n$.
		\item
		If $S$, $\tilde S$ satisfy condition (U), then $(\bb x_0 + S) \cap (\bb x^* \circ \exp \tilde S^\perp)$ contains {\bf at most one} point, for any $\bb x_0$, $\bb x^* \in \rrpp^n$.
		\end{enumerate}
	\end{prb}

If a generalized mass-action system happens to be a classical mass-action system, then its stoichiometric subspace $S$ is also the kinetic-order subspace $\tilde S$. The existence and uniqueness of a point in the intersection $(\bb x_0 + S) \cap (\bb x^* \circ \exp S^\perp)$ for any $\bb x_0$, $\bb x^* \in \rrpp^n$ is the content of Birch's theorem in algebraic statistics~\cite{SturmfelsPachter2005}.
%%%%%%%%%%%%%%%%%%%%%%%%%%%%%%

\bigskip

	Another reformulation of the above problems was introduced by M\"uller and Regensburger~\cite{MR12}, in terms of injectivity/surjectivity of an exponential map or a generalized polynomial map onto a polyhedral cone. Such polynomial maps appear in other applications; for example, a renormalized version of the generalized polynomial appears in computer graphics and geometric modelling, where the map being injective implies that the curve or surface does not self-intersect~\cite{geommodel}.

	Let $\bb x^* \in \rrpp^n$ be an arbitrary vector, and $S$, $\tilde S \subseteq \rr^n$ be vector subspaces, with $d =\codim S$, $\tilde d = \codim \tilde S$. Choose a basis for $S^\perp$ and let the basis vectors be the rows of the matrix $W \in \rr^{d\times n}$. Similarly, choose a basis for $\tilde S^\perp$ and let the basis vectors be the rows of $\tilde W \in \rr^{\tilde d \times n}$. Write the two matrices in terms of their columns: $W = (\bb w^1, \bb w^2, \cdots, \bb w^n)$ and $W = (\tilde{\bb w}^1, \tilde{\bb w}^2,\cdots, \tilde{\bb w}^n)$. In this manner, $S^\perp = \Img (W^T)$, $S = \Ker W$, and $\tilde S^\perp = \Img (\tilde W^T)$, $\tilde S = \Ker \tilde W$. Finally, write $\intr{C}(W)$ for the relative interior of the polyhedral cone $C(W)$, i.e., $\intr{C}(W)$ is the set of all positive combinations of $\{\bb w^i \}_{i=1}^n$. For any $\bb x^* \in \rrpp^n$, define the maps
	\eq{
		\begin{array}{r}
			f_{\bb x^*} :
			\\ \hphantom{h}
		\end{array} &
		\begin{array}{ccccc}
		\rrpp^{\tilde d} & \to & \intr{C}(W) \subseteq \rr^d , \\
			 \bb \xi & \mapsto & W(\bb x^* \circ {\bb \xi}^{\tilde W})
			& = &\sum_{i=1}^n x^*_i \xi^{\tilde{\bb w}^i} \bb w^i,
		\end{array}
		\intertext{and}
		\begin{array}{r}
			F_{\bb x^*} :
			\\ \hphantom{h}
		\end{array}&
		\begin{array}{ccccc}
		\rr^{\tilde d} & \to & \intr{C}(W) \subseteq \rr^d ,\\
			 \bb \lambda & \mapsto & W(\bb x^* \circ e^{\tilde W^T{\bb \lambda}}) & = & \sum_{i=1}^n x^*_i e^{\braket{\tilde{\bb w}^i}{\lambda}} \bb w^i.
		\end{array}
	}
Problem~\ref{prb:mfld} is equivalent to the following (see \cite{MR12, MR14} for details):
	\begin{prb}[3]
	\label{prb:bij}
		Let $S$, $\tilde S \subseteq \rr^n$ be vector subspaces. %For any vector $\bb x^* \in \rrpp^n$, define matrices $W$, $\tilde W$ and function $f_{\bb x^*}$ (respectively $F_{\bb x^*}$) as above.
What are conditions (E) and (U) on $S$, $\tilde S$, so that the following statements are true?
		\begin{enumerate}\itemsep=-2pt
		\item
		If $S$, $\tilde S$ satisfy condition (E), then the map $f_{\bb x^*}$ (respectively $F_{\bb x^*}$) is {\bf surjective} onto $C^0(W)$, for any $\bb x^* \in \rrpp^n$.
		\item
		If $S$, $\tilde S$ satisfy condition (U), then the map $f_{\bb x^*}$ (respectively $F_{\bb x^*}$) is {\bf injective}, for any $\bb x^*\in \rrpp^n$.
		\end{enumerate}
	\end{prb}

\medskip

	 M\"uller and Regensburger characterized when the maps $f_{\bb x^*}$, $F_{\bb x^*}$ are injective, namely, if and only if $\sigma(S) \cap \sigma(\tilde S^\perp) = \{\bb 0\}$~\cite[Theorem 3.6]{MR12}. Recall that, for a subset $S \subseteq \rr^n$, its set of \emph{sign vectors} $\sigma(S)$ is the image of vectors in $S$ under the coordinate-wise sign function (Definition~\ref{def:signvectors}). They also provided a \emph{sufficient} condition for bijectivity: if $\sigma(S) = \sigma(\tilde S)$ and $(+,+,\cdots,+)^T \in \sigma(S^\perp)$, then $f_{\bb x^*}$, $F_{\bb x^*}$ are bijective (and indeed, real analytic isomorphisms)~\cite[Proposition 3.9]{MR12}. Our main result (Theorem~\ref{thm:mainCRN}) can be regarded as a generalization of this result. Recently, M\"uller, Hofbauer, and Regensburger have used Hadamard's global inversion theorem to \emph{characterize} when $f_{\bb x^*}$, $F_{\bb x^*}$ are bijective for arbitrary $\bb x^* \in \rrpp^n$~\cite{HMR}.

%	\begin{thm}[{\cite[Theorem ??]{addref17}}]
%	\label{thm:addref17}
%	Let $S, \tilde S \subseteq \rr^n$ be vector subspaces. For any vector $\bb x^* \in \rrpp^n$, define matrices $W, \tilde W$ and functions $f_{\bb x^*}, F_{\bb x^*}$ as above. Then $f_{\bb x^*}$ (respectively $F_{\bb x^*}$) is bijective if and only if
%		\begin{enumerate}\itemsep=-2pt
%		\item
%			$\sigma(S) \cap \sigma(\tilde S^\perp) = \{\bb 0\}$;
%		\item
%			for any non-zero sign vector $\tilde \tau \in \sigma(\tilde S^\perp)_{\geq 0}$, there is a non-zero sign vector $\tau \in \sigma(S^\perp)_{\geq 0}$ such that $\tau \leq \tilde \tau$, and
%		\item
%			the pair $(S, \tilde S)$ is {\emph{non-degenerate}}.\footnote{This non-degeneracy condition is a condition on set of sign vectors $\sigma(S)$ and the vector space $\tilde S$.}
%		\end{enumerate}
%	\end{thm}

%%%%%%%%%%%%%%%%%%%%%%%%%%%%%%%%%%%%%%%%%%%%%%%%%%%%%%%%%%%%%%%%%%%%%%%%%%%%%%%%%%%%%%%%%%%%%%%%%%%%%%%%%%%%%%%%%%%%%%%
\section{Main result}
\label{sec:Result}

	In previous work as well as in ours, the conditions (E) and (U) are stated in terms of {{sign vectors}}. For a brief introduction to sign vectors of linear subspaces, we refer the reader to the appendix in \cite{MR12}.

\begin{defn}
\label{def:signvectors}
		Given a vector $\bb x \in \rr^n$, we define its {\df{sign vector}} to be
		\eqn{
			\sigma(\bb x) =
			(\sgn(x_1) , \sgn(x_2) , \cdots , \sgn(x_n))^T
			\in \{0,+,-\}^n.
		}
The set of sign vectors for a subset $S \subseteq \rr^n$ is the collection $\sigma(S) = \{ \sigma(\bb x) : \bb x \in S\}$.

We introduce a partial order on $\{0,+,-\}$ by $0<-$ and $0<+$ (but no relation between $-$ and $+$). This induces a partial order on $\{0,+,-\}^n$: $\tau \leq \tau'$ if $\tau_j \leq \tau'_j$ for all $j$. The {\df{closure}} of a set of sign vectors $\Lambda \subseteq \{0,+,-\}^n $ is the set
		\eqn{
			\cl{\Lambda} = \{ \tau : \text{ there exists } \tau' \in \Lambda \text{ such that } \tau \leq \tau' \}.
		}
\end{defn}

\bigskip

	We define an \emph{orthant}\footnote{This differs from the typical definition of an orthant of $\rr^n$, which is full dimensional.
	} %
 of $\rr^n$ to be a maximal subset of $\rr^n$ on which $\sigma$ is constant. Geometrically, the sign vector $\sigma(\bb x)$ tells us which orthant $\mc O_{\bb x}$ the vector $\bb x$ lies in, while the closure $\cl{\sigma(\bb x)}$ refers to the union of $\mc O_{\bb x}$ and its boundary.  Finally, we define an orthogonality relation on $\{0, + , - \}^n$; we say that two sign vectors $\tau$ and $\tau'$ are {\df{orthogonal}} (denoted $\tau \perp \tau'$) if
	\begin{enumerate}[label=\roman*.]\itemsep=-3pt
	\item[]
		either $\tau_j\cdot \tau_j' = 0$ for all $1 \leq j \leq n$
	\item[]
		or there exist indices $i, j$ such that $\tau_i\cdot \tau_i' = +$ and $\tau_j\cdot \tau_j' = -$,
	\end{enumerate}
where the product operation on signs is as one would expect:
	\eq{
		 + \cdot + = - \cdot - = +, \qquad
		+ \cdot - = -, \quad \text{and} \quad
		 + \cdot\, 0 = - \cdot 0 = 0 \cdot 0 = 0.
	}
It is easy to see that if $\bb x$, $\bb y \in \rr^n$ are orthogonal vectors, then $\sigma(\bb x) \perp \sigma(\bb y)$. \\

\bigskip

	We show in this section that if $\sigma(S) \subseteq \cl{\sigma(\tilde S)}$ and $\dim S = \dim \tilde S$, then for any $\bb x_0$, $\bb x^* \in \rrpp^n$, the intersection $(\bb x_0 +S) \cap (\bb x^*\circ \exp \tilde S^\perp)$ contains exactly one point. The intuitive idea is that the sign condition $\sigma(S) \subseteq \cl{\sigma(\tilde S)}$ is related to a transversal intersection of the two manifolds $(\bb x_0 + S)$ and $(\bb x^*\circ \exp \tilde S^\perp)$. If we have one intersection point, say $\bb x^* \in (\bb x^*+S) \cap (\bb x^*\circ \exp\tilde S^\perp)$, we cannot lose the intersection point as we translate the affine plane from $(\bb x^* +S)$ to $(\bb x_0 +S)$.

	We first show in Lemma~\ref{lem:uniqueness} that our sign condition $\sigma(S) \subseteq \cl{\sigma(\tilde S)}$ implies the uniqueness condition $\sigma(S) \cap \sigma(\tilde S^\perp) = \{ \bb 0\}$ in \cite{MR12}. %\cite[Proposition 3.1]{MR12}.
 In Lemma~\ref{lem:transversal}, we establish transversality of the manifolds $(\bb x + S)$ and $(\bb x^*\circ \exp \tilde S^\perp)$. Lemma~\ref{lem:compactness} prevents our desired intersection point from escaping to the boundary of $\rrpp^n$ or to infinity. Finally, these results lead to Theorem~\ref{thm:main}, concluding the existence and uniqueness of a point in the intersection $(\bb x_0 + S) \cap (\bb x^* \circ \exp\tilde S^\perp)$. In Theorem~\ref{thm:mainCRN} and Corollary~\ref{cor:last}, we apply this result to generalized mass-action systems.

\begin{lem}[Uniqueness]
\label{lem:uniqueness}
	Let $S$, $\tilde S \subseteq \rr^n$ be vector subspaces. If $\sigma(S) \subseteq \cl{\sigma(\tilde S)}$, then $\sigma(S) \cap \sigma(\tilde S^\perp) = \{\bb 0\}$. In particular, for any $\bb x_0$, $\bb x^* \in \rrpp^n$ the intersection $(\bb x_0 + S) \cap (\bb x^* \circ \exp \tilde S^\perp)$ contains at most one point.
\end{lem}

	\proof
	By assumption, $\sigma(S) \cap \sigma(\tilde S^\perp) \subseteq \cl{\sigma(\tilde S)} \cap \sigma(\tilde S^\perp)$. We show that $\cl{\sigma(\tilde S)} \cap \sigma(\tilde S^\perp) = \{ \bb 0 \}$. Let $\tau \in \cl{\sigma(\tilde S)} \cap \sigma(\tilde S^\perp)$ be a sign vector. There exist vectors $\bb x \in \tilde S$  and $\bb y \in \tilde S^\perp$ such that $\tau \leq \sigma(\bb x)$ and $\tau = \sigma(\bb y)$. It is easy to see that if $\tau \leq \sigma(\bb x)$, and $\tau \perp \sigma(\bb x)$, then $\tau = \bb 0$.

%	In general, note that if $\tau' \leq \tau''$ and $\tau'$ and $\tau''$ are orthogonal, then $\tau' = 0$. Hence $\tau = 0$.

%	Since $\bb x$ and $\bb y$ are orthogonal, their sign vectors are orthogonal. For each of its non-zero coordinates $\tau_j \neq 0$, we know that $\sigma(\bb x)_j = \tau_j$; thus there cannot be an index such that $\tau_j \cdot \sigma(\bb x)_j = -$. Hence, $0 = \tau_j \cdot  \sigma(\bb x)_j = \tau_j\cdot  \tau_j$ for $1 \leq j \leq n$, and $\tau = \bb 0$.

	By \cite{MR12}, $\sigma(S) \cap \sigma(\tilde S^\perp) = \{\bb 0\}$ is necessary and sufficient for the intersection $(\bb x_0+S)\cap (\bb x^*\circ \exp \tilde S^\perp)$ to contain at most one point for any $\bb x_0$, $\bb x^* \in \rrpp^n$.
	\qed \\

\begin{lem}[Compactness]
\label{lem:compactness}
	Let $S$, $\tilde S \subseteq \rr^n$ be vector subspaces, and let $K \subseteq \rrpp^n$ be a compact subset, and $\bb x^* \in \rrpp^n$. Suppose $\sigma(S) \subseteq \cl{\sigma(\tilde S)}$. Then $(K+S) \cap (\bb x^* \circ \exp \tilde S^\perp)$ is a compact subset of $\rrpp^n$.
\end{lem}

	\proof
	Let $\Gamma = (K+S) \cap (\bb x^* \circ \exp \tilde S^\perp)$. Since $\bb x^* \circ \exp \tilde S^\perp \subseteq \rrpp^n$, the intersection $\Gamma$ also lies in the positive orthant. We first show that $\Gamma$ is bounded away from infinity and from the boundary of $\rrpp^n$.

	Suppose that is not the case. Let $\bb x^k \in \Gamma$ be a sequence such that either $ \limsup_{k\to\infty} x^k_i = \infty$ or $ \liminf_{k\to\infty} x^k_i = 0$ for some index $1 \leq i \leq n$. Passing to a subsequence, we may assume that
	\eq{
		\begin{array}{lc}
		\ds \lim_{k\to\infty} x^k_i = \infty  & \qquad \text{for } i \in I_1, \\
		\ds \lim_{k\to\infty} x^k_i = 0 &  \qquad \text{for } i \in I_2, \\
		\ds \lim_{k\to\infty} x^k_i \in (0,\infty) & \qquad \text{for } i \in I_3,
		\end{array}
	}
where $I_1$, $I_2$, $I_3$ partition the index set $\{1,2,\dots, n\}$, and $I_1 \union I_2 \neq \emptyset$.

	On one hand, $\bb x^k \in K+S$, so decompose it as $\bb x^k = \bb v^k + \bb s^k$, where $\bb v^k \in K$ and $\bb s^k \in S$. Since $K \subseteq \rrpp^n$ is compact, each component of $\bb v^k$ is uniformly bounded from above and below from zero. Thus for $i \in I_1$, we have $s^k_i \to \infty$; in particular, $s^k_i > 0$ for sufficiently large $k$. Similarly, if $i \in I_2$, then $s^k_i < 0$ for sufficiently large $k$, because $s^k_i + v^k_i \to 0$ and $v^k_i > 0$ is bounded away from zero. Hence the sign of $s^k_i$ is constant for any $i \in I_1 \union I_2$ for sufficiently large $k$. Because $\sigma(\bb s^k) \in \sigma(S) \subseteq \cl{\sigma(\tilde S)}$, there is a vector $\bb u \in \tilde S$ such that $u_i > 0$ if $i \in I_1$ and $u_i < 0$ if $i \in I_2$.

	On the other hand, $\bb x^k \in \bb x^* \circ \exp \tilde S^\perp$, that is, $\log \left( \frac{\bb x^k}{\bb x^*}\right) \in \tilde S^\perp$, where the division is understood to be component-wise. Hence, $\bb u \perp \log \left( \frac{\bb x^k}{\bb x^*} \right)$ for all $k$, and we have
	\eq{
		0 &=
		\left<\bb u \,,\, \log \left( \frac{\bb x^k}{\bb x^*}\right) \right>
		= \sum_{i\in I_1} u_i \log \left( \frac{x^k_i}{x^*_i} \right)
			+ \sum_{i\in I_2} u_i \log \left( \frac{x^k_i}{x^*_i} \right)
			+ \sum_{i\in I_3} u_i \log \left( \frac{x^k_i}{x^*_i} \right).
	}
The sum over $I_3$ is uniformly bounded for all $k$. Now let $k\to \infty$. For $i \in I_1$, we know $u_i > 0$ and $x_i^k \to \infty$, so the sum over $I_1$ is positive and unbounded. For $i \in I_2$, we know $u_i < 0$ and $x_i^k \to 0$, so $\log\left(\frac{x_i^k}{x_i^*} \right) \to -\infty$, so the sum over $I_2$ is also positive and unbounded. Consequently, $0 = \lim_{k\to\infty} \braket{\bb u}{\log\left( \frac{\bb x^k}{\bb x^*} \right)} = \infty$, a contradiction. Hence, $\Gamma \subseteq \rrpp^n$ is bounded away from infinity and away from the boundary of the positive orthant.

	Next, we want to show that $\Gamma \subseteq \rrpp^n$ is a closed subset. Let us fix $\epsilon > 0$ such that $\Gamma$ lies inside the hypercube $Q = [\epsilon, \epsilon^{-1}]^n \subseteq \rrpp^n$. Being the intersection of two closed sets, $Q \cap (K+S)$ is closed. The set $Q \cap (\bb x^* \circ \exp \tilde S^\perp)$ is diffeomorphic to $[\log\epsilon, \log\epsilon^{-1}]^n \cap (\ln \bb x^* + \tilde S^\perp)$, which is again a closed set. Therefore, the set $(K+S)\cap (\bb x^* \circ \exp \tilde S^\perp) = [Q \cap (K+S)] \cap [Q \cap (\bb x^* \circ \exp \tilde S^\perp) ]$ is the intersection of two closed sets, and thus it is closed in $\rrpp^n$.
	\qed \\

	 Two manifolds $X$ and $Y$ of $\rr^n$ {\df{intersect transversally}} if at each point $\bb p \in X \cap Y$, their tangent spaces span the entire Euclidean space, i.e., $T_{\bb p}(X) + T_{\bb p}(Y) = \rr^n$. We refer the reader to \cite{difftop, difftop2} for the theory of transversality and intersection.

	Again, let $\bb x_0$, $\bb x^*\in \rrpp^n$ be two arbitrary vectors in what follows. In Lemma~\ref{lem:uniqueness}, we showed that our sign condition $\sigma(S) \subseteq \cl{\sigma(\tilde S)}$ implies $\sigma(S) \cap \sigma(\tilde S^\perp) = \{\bb 0\}$, which is equivalent to the intersection $(\bb x_0+S) \cap (\bb x^*\circ \exp \tilde S^\perp)$ containing at most one point. Indeed, this weaker sign condition together with $ \dim S = \dim \tilde S $ is enough to conclude that the two manifolds $\bb x_0+S$ and $\bb x^*\circ \exp \tilde S^\perp$ intersect transversally. This is the content of the follow lemma.

\begin{lem}[Transversality]
\label{lem:transversal}
	Let $S$, $\tilde S \subseteq \rr^n$ be vector subspaces. Assume $\sigma(S) \cap \sigma(\tilde S^\perp) = \{ \bb 0\}$. Let $\bb x_0$, $\bb x^* \in \rrpp^n$ be any two positive vectors. Then the tangent spaces of $\bb x_0 +S$ and $\bb x^*\circ \exp \tilde S^\perp$ satisfy
	\eq{
		T_\bb p(\bb x_0 +S) \cap T_\bb p(\bb x^*\circ \exp \tilde S^\perp) = \{\bb 0\}
	}
for any point $\bb p \in (\bb x_0 +S) \cap (\bb x^*\circ \exp \tilde S^\perp)$.

If we further assume that $\dim S = \dim \tilde S$, then $T_\bb p(\bb x_0 + S) + T_\bb p(\bb x^*\circ \exp \tilde S^\perp) = \rr^n$ for any intersection point $\bb p \in (\bb x_0+S) \cap (\bb x^*\circ \exp \tilde S^\perp)$, i.e., $\bb x_0 +S$ and $\bb x^* \circ \exp \tilde S^\perp$ intersect transversally.
\end{lem}

	\proof
	For any intersection point $\bb p \in (\bb x_0 +S) \cap (\bb x^*\circ \exp \tilde S^\perp)$, we note that $T_\bb p(\bb x_0 + S) = S$ and $T_\bb p(\bb x^*\circ \exp \tilde S^\perp) = \bb p \circ \tilde S^\perp$ and hence $\sigma(T_\bb p(\bb x^*\circ \exp \tilde S^\perp)) = \sigma(\tilde S^\perp)$. 
	
	Now consider $\bb x \in T_\bb p(\bb x_0 + S) \cap T_\bb p(\bb x^*\circ \exp \tilde S^\perp)$. Then $\sigma(\bb x) \in \sigma(S) \cap \sigma(\tilde S^\perp)  = \{\bb 0\}$, which implies $\bb x = \bb 0$. Consequently, $T_\bb p (\bb x_0 + S) \cap T_\bb p(\bb x^* \circ \exp \tilde S^\perp) = \{ \bb 0\}$.

	If we further assume that $\dim S = \dim \tilde S$, we note that $T_\bb p (\bb x_0 +S ) + T_\bb p (\bb x^*\circ \exp \tilde S^\perp)$ is of dimension  %$d + (n-d) = n$
$n$. In other words, the manifolds $\bb x_0+S$ and $\bb x^*\circ \exp \tilde S^\perp$ intersect transversally.
	\qed \\

%%%%%%%%%%%%%%%%%%%%%%%%%%%%%%

	Now we are ready to state and prove our main result. The proof starts with a known intersection point, $\bb x^* \in (\bb x^* +S)\cap (\bb x^* \circ \exp \tilde S^\perp)$. Next, we translate the affine space $(\bb x^* + S)$ to $(\bb x_0 +S)$, creating a $(d+1)$-dimensional strip of the form $K+S$, where $d = \dim S$ and $K$ is a compact subset of $\rrpp^n$. This strip  intersects $\bb x^* \circ \exp \tilde S^\perp$ transversally, and we use Corollary~\ref{thm:difftop2} below to conclude that the intersection $(K+S) \cap (\bb x^* \circ \exp \tilde S^\perp)$ is a one-dimensional manifold, whose boundary lies on the boundary of the affine strip $K+S$. Finally, we conclude the existence of a boundary point on $\bb x_0 +S$ by the uniqueness condition. \\

%%%{\orange{%%%---%%%
%%% We will use 
%%%}}%---%%%
%%%the following differential topology result:
%%%
%%%	\begin{thm}[\cite{difftop}] %page 60
%%%	\label{thm:difftop}
%%%		Let $X$ be a submanifold of $Y$ with boundary, and let $Z$ be a boundaryless submanifold of $Y$. Suppose $X$ intersects $Z$ transversally, and $\bdy X$ intersects $Z$ transversally. Then $X \cap Z$ is a manifold with boundary $\bdy(X\cap Z) = Z \cap \bdy X$, and $\codim_X(X\cap Z) = \codim_Y(Z)$.
%%%	\end{thm}
%%%
%%%
%%%	For our purpose, the ambient manifold is $Y = \rrpp^n$, while the submanifold with boundary is $X = K+S$ for some compact set $K$ and the boundaryless submanifold is $Z = \bb x^* \circ \exp \tilde S^\perp$.
%%%
%%%%%%---%%%
%%%\begin{tikzpicture}[overlay]
%%%	\node at (0,0) {};
%%%	\draw [red, very thick] (-1,0)--(-1,3.3) node [midway, left] {replace with};
%%%	\node at (-1, 1.3) [left] {\red{ green text }};
%%%	\node at (-1,1.) [left] {\red{below?}};
%%%\end{tikzpicture}

	Consider the following differential topology result:
	\begin{thm}[{\cite[Theorem~3.5.1]{difftop}}]
	\label{thm:difftop}
		Let $X$ and $Y$ be manifolds and $Z \subseteq Y$ a submanifold, where $Z$ and $Y$ are boundaryless. Let $f: X \to Y$ be a smooth map. Suppose $f$ intersects $Z$ transversally and $\left. f \right|_{\bdy X}$ also intersects $Z$ transversally. Then $f^{-1}(Z)$ is a submanifold of $X$ with boundary $\bdy(f^{-1}(Z)) = \bdy X \cap f^{-1}(Z)$, and $\codim_{X}(f^{-1}(Z)) = \codim_Y(Z)$. 
	\end{thm}
	
%For our purpose, the ambient manifold is $Y = \rr^n_{>0}$, while $X =  K+S$, for some compact set $K \subseteq \rr^n_{>0}$ is a manifold with boundary $\bdy(X) = (\bb x_0 +S) \cup (\bb x^* +S)$, and $Z = \bb x^* \circ \exp \tilde S^\perp$ is a boundaryless submanifold of $Y$. If $f$ is the inclusion map of $K+S$ into $\rr^n_{>0}$, to say that the maps $f$ and $\left. f\right|_{\bdy X}$ intersect the manifold $Z$ transversally is equivalent to the manifolds $X$ and $\bdy X$ intersect $Z$ transversally. Therefore, the version that we will apply directly is the following:

Consider the setting where the ambient manifold is $Y = \rr^n_{>0}$. If $f$ is the inclusion map of a submanifold $X$ into $\rr^n_{>0}$, to say that the maps $f$ and $\left. f\right|_{\bdy X}$ intersect the manifold $Z$ transversally is equivalent to the manifolds $X$ and $\bdy X$ intersect $Z$ transversally. The preimage $f^{-1}(Z)$ is the submanifold $X \cap Z$. Moreover, the dimension of the intersection $X \cap Z$ is given by the equation
	\eq{
		\dim X - \dim (X \cap Z) = \codim_X(X\cap Z) = \codim_{\rr^n_{>0}}(Z) = n - \dim Z.
	}
In other words, $\dim(X\cap Z) = \dim X + \dim Z - n$. We arrive at the following corollary:

	\begin{cor}
	\label{thm:difftop2}
		Let $X$, $Z \subseteq \rr^n_{>0}$ be submanifolds, where $Z$ is boundaryless. Suppose $X$ intersects $Z$ transversally and $\bdy X$ also intersects $Z$ transversally. Then $X \cap Z$ is a manifold with boundary $\bdy(X \cap Z) = \bdy X \cap Z$ and of dimension $\dim (X \cap Z) = \dim X + \dim Z - n$.
	\end{cor}

%%
%%%${}$\\
%%\noindent
%%{\red{%%%---%%%
%%{\bf Note to ourselves:} \\
%%-- The theorem in Guillemin and Pollack~\cite{difftop} at the bottom of page 60 also uses the language of $f$ and $Z$ intersecting transversally. Their explanation of using $f = \iota$ the inclusion map is at the bottom of page 29.\\
%%${}$\\
%%Second important note is that their theorem says $f: X \to Y$ is an onto map, which is not true in our case. (We can drop the assumption and the theorem is still holds, but that is not what is printed in the book.)\\
%%${}$\\
%%-- A different reference: Elements of Differential Topology, by Anant R. Shastri.  Here, again theorem uses the language of $f$ and $Z$ intersecting transversally. This is Theorem 3.5.1 (Transversal Inverse Image Theorem) on page 94. Their connecting with using $f = \iota$ is in Exercise 3.5 question 2 on page 95. \\
%%${}$\\
%%-- The Guillemin-Pollack reference has the advantage that the book does emphasize geometric intuition a lot, and seems not too difficult to read. The Anant book has the advantage that their theorems are numbered and also not having the ``onto'' assumption. 
%%}}%---%%%

\bigskip

Our main result is:

\begin{thm}
\label{thm:main}
	Let $S$, $\tilde S \subseteq \rr^n$ be vector subspaces of equal dimension with $\sigma(S) \subseteq \cl{\sigma(\tilde S)}$. Then for any positive vectors $\bb x_0$, $\bb x^* \in \rrpp^n$, the intersection $(\bb x_0+S) \cap (\bb x^*\circ \exp \tilde S^\perp)$ consists of exactly one point.
\end{thm}

	\proof
		Let $\bb x_0$, $\bb x^* \in \rrpp^n$ be arbitrary positive vectors. Lemma~\ref{lem:uniqueness} implies that the intersection $(\bb x_0+S) \cap (\bb x^*\circ \exp \tilde S^\perp)$ contains at most one point. Consider first $\bb x^* \in \bb x_0 +S$. Clearly, $(\bb x_0 +S) \cap (\bb x^* \circ \exp \tilde S^\perp) = \{ \bb x^*\}$.

	Now consider the case when $\bb x^* \not\in \bb x_0 +S$. Let $d = \dim S$. We define a $(d+1)$-dimensional affine strip, which we will intersect with $(\bb x^* \circ \exp \tilde S^\perp)$. To define this affine strip, we consider the interpolation function
		\eq{
			\begin{array}{cccl}
				K: & [0,1] & \to & \rrpp^n, \\
					& \delta & \mapsto  &\delta \bb x_0 + (1-\delta) \bb x^*.
			\end{array}
		}
Since the line segment $K([0,1]) \subseteq \rrpp^n$ is compact, the intersection $(K([0,1]) + S) \cap (\bb x^*\circ \exp \tilde S^\perp) \subseteq \rrpp^n$ is compact by Lemma~\ref{lem:compactness}. Moreover, the manifolds $K([0,1])+S$ and $\bb x^*\circ \exp \tilde S^\perp$ intersect transversally, as a consequence of Lemma~\ref{lem:transversal}, i.e.,
	\eq{
		T_\bb p(K([0,1])+S) + T_\bb p(\bb x^*\circ \exp \tilde S^\perp)
		\supseteq T_\bb p(\bb x^*+S) + T_\bb p(\bb x^*\circ \exp \tilde S^\perp)
		= \rr^n.
	}
By Corollary~\ref{thm:difftop2}, the intersection $\Gamma = (K([0,1])+S)\cap (\bb x^* \circ \exp\tilde S^\perp)$ is a manifold with boundary $\bdy \Gamma \subseteq \bdy (K([0,1])+S) = (\bb x^* + S) \cup (\bb x_0 +S)$. In addition, $\Gamma$ is 1-dimensional because 
	\eq{
		\dim(\Gamma) = \dim(K([0,1])+S) + \dim(\bb x^* \circ \exp \tilde S^\perp) - n 
			= 1 + \dim S + \dim \tilde S^\perp - n = 1.
	} 

		Consider the connected component $\Gamma^* \subseteq \Gamma$ containing the point $\bb x^*$. The point $\bb x^*$ must be an endpoint of $\Gamma^*$; otherwise uniqueness fails at $K(\delta_0) +S$ for some small $\delta_0 > 0$. Since $\Gamma^*$ is compact, it is a curve with two endpoints. As $\bdy \Gamma^* \subseteq \bdy \Gamma = (\bb x^* + S) \cup (\bb x_0 + S)$, by uniqueness the other endpoint of $\Gamma^*$ must be in $\bb x_0 +S $. Thus, a point exists in $(\bb x_0 +S) \cap (\bb x^* \circ \exp \tilde S^\perp)$. 
	\qed \\

We apply Theorem~\ref{thm:main} to show the existence and uniqueness of vertex-balanced steady state for a generalized mass-action system.

\begin{thm}[Vertex-balanced steady states of a generalized mass-action system]
\label{thm:mainCRN}
	Let $(G, \Phi, \tilde \Phi)$ be a 
	weakly reversible
	generalized reaction network, with stoichiometric subspace $S$ %, and suppose that every vertex of $G$ is the source of some edge, so that the 
	and
	kinetic-order subspace $\tilde S$. %is well-defined. 
	Assume that $\dim S = \dim \tilde S$ and $\sigma(S) \subseteq \cl{\sigma(\tilde S)}$. Then the following statements hold:
	\begin{enumerate}[label=\roman*)]
	\item
		Suppose for some rate constants ${\boldsymbol\kk}$, the generalized mass-action system $(G_{\boldsymbol\kk}, \Phi, \tilde \Phi)$ admits a vertex-balanced steady state $\bb x^*$. Then every stoichiometric compatibility class contains exactly one vertex-balanced steady state.
	\item
		%Suppose $G$ is weakly reversible. Then  
		$\tilde \delta_G = 0$ if and only if the generalized mass-action system ($G_{\boldsymbol\kk}, \Phi, \tilde \Phi)$ admits a vertex-balanced steady state $\bb x^*$ for all rate constants ${\boldsymbol\kk}$. %Moreover, 
		In this case,
		every stoichiometric compatibility class contains exactly one vertex-balanced steady state.
	\item
		Under the premises of i), additionally suppose $\delta_G = 0$. Then every stoichiometric compatibility class contains exactly one positive steady state, which is vertex-balanced.
	\end{enumerate}
\end{thm}

	\proof
		As $\bb x^*$ is a vertex-balanced steady state for $(G_{\boldsymbol \kk}, \Phi, \tilde \Phi)$, the set of vertex-balanced steady state is $Z_{\boldsymbol \kk} = \bb x^* \circ \exp \tilde S^\perp$. By Theorem~\ref{thm:main}, $Z_{\boldsymbol \kk}$ intersects the stoichiometric compatibility class $\bb x_0 +S$ exactly once for any $\bb x_0 \in \rrpp^n$. This proves statement {\it i)}.

        The first part of statement {\it ii)} is the content of \cite[Theorem~1(a)]{MR14}.
	%If $G$ is weakly reversible and $\tilde \delta_G = 0$, then the set of vertex-balanced steady states $Z_{\boldsymbol \kk} \neq \emptyset$ for any $\boldsymbol \kk > 0$~\cite{MR14}. 
	By statement {\it i)}, we conclude that every stoichiometric compatibility class contains exactly one vertex-balanced steady state. %The converse of (ii) was also proved in~\cite{MR14}.

		If in addition, $\delta_G = 0$, then $E_{\boldsymbol\kk} = Z_{\boldsymbol\kk}$, i.e., there are no positive steady states that are not vertex-balanced. Consequently, there exists a unique steady state within each stoichiometric compatibility class, which is vertex-balanced. This proves statement {\it iii)}.
	\qed \\

We state a simpler version of {\it iii)} in the theorem above.

\begin{cor}
\label{cor:last}
	Let $(G,\Phi,\tilde \Phi)$ be a weakly reversible generalized reaction network, with stoichiometric subspace $S$ and kinetic-order subspace $\tilde S$. Suppose that $\dim S =\dim \tilde S$, $\sigma(S) \subseteq \cl{\sigma(\tilde S)}$, and $\delta_G = \tilde \delta_G = 0$. Then for any choice of rate constants, every stoichiometric compatibility class contains exactly one positive steady state, which is vertex-balanced.
\end{cor}

%\rmk
%	Our condition is neither more nor less general than the sufficient condition for existence of intersection point provided by M\"uller and Regensburger in \cite{MR12}. Their result does not apply to the example in Sec.~\ref{sec:MainEx}, while our result is not relevant for networks that admit multiple vertex-balanced steady states. For example, consider the case when
%	\eq{
%		S &= \Span_\rr \left\{
%			\begin{pmatrix} -2 \\ -1 \\ 2 \\ 0 \end{pmatrix},
%			\begin{pmatrix} -1 \\ -2 \\ 0 \\2 \end{pmatrix}
%			\right\}, \qquad
%		%\\
%		\tilde S = \Span_\rr \left\{
%			\begin{pmatrix} -1 \\ -2 \\ 2 \\ 0 \end{pmatrix},
%			\begin{pmatrix} -2 \\ -1 \\ 0 \\2 \end{pmatrix}
%			\right\}.
%	}

%\section{Multistability}
%\label{sec:multistability}
%
%	We conclude by
%
%
%	\begin{prop}
%	\label{prop:multistability}
%		Let $(G_\kk, \Phi, \tilde \Phi)$ be a weakly reversible generalized mass-action system. Assume that the stoichiometric subspace $S$ and kinetic-order subspace $\tilde S$ satisfy the sign condition $\sigma(S) \cap \sigma(\tilde S^\perp) = \{\bb 0\}$. If $\delta = 0$, the system is not multistable.
%	\end{prop}
%
%
%	-- Note that $\sigma(S) \subseteq \cl{\sigma(\tilde S)}$ implies existence of unique vertex-balanced steady state, but the converse is not true.
%
%	-- more subtlety: we don't know that in general $\tilde Z_\kk = \tilde E_\kk$, unless $\delta = 0$.

	\bigskip

We have focused almost exclusively on the {\em existence} and uniqueness of vertex-balanced steady states for generalized mass-action systems. For complex-balanced equilibria of classical mass-action systems, more is known. For example, complex-balanced equilibria are locally asymptotically stable within their stoichiometric compatibility classes. They are conjectured to be globally stable in their stoichiometric compatibility classes; this is known as the \emph{global attractor conjecture}~\cite{CraciunGAC, Toric09}. In particular, it has been shown that a complex-balanced equilibrium of a mass-action system is globally stable within its stoichiometric compatibility class if the network has a single connected component~\cite{DaveGAC}, or is strongly endotactic~\cite{CasianGAC1, CasianGAC2}, or if the system is in $\rr^3$~\cite{CasianGAC1, CasianGAC2}. A proof of the global attractor conjecture in full generality has been proposed in \cite{CraciunGAC}.

For {\em planar} generalized mass-action systems (in particular, S-systems),
local and even global stability of vertex-balanced steady states have been characterized in \cite{BorosHofbauerMueller2017,BorosHofbauerMuellerRegensburger2018,BorosHofbauerMuellerRegensburger2019}.
For generalized mass-action systems of arbitrary dimension, necessary conditions for linear stability have been given in \cite{BorosMuellerRegensburger2019}. Obviously, it is not true that a vertex-balanced steady state is always globally stable  within its stoichiometric compatibility class, since it is possible for a generalized mass-action system to have multiple vertex-balanced steady states within the same stoichiometric compatibility class. Consider, for example, the following generalized mass-action system:
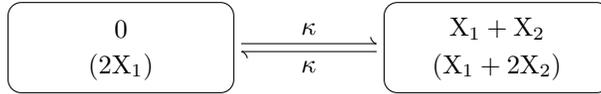
\tikzc{
		\draw [rounded corners=5pt] (0,0) rectangle++(3,1.2);
		\draw [rounded corners=5pt] (5,0) rectangle++(3,1.2);
		\node at (1.5,0.85) {$0$};
		\node at (1.5,0.35) {($\chemfig{2X_1}$)};
		\node at (6.5,0.85) {$\chemfig{X_1 + X_2}$};
		\node at (6.5,0.35) {($\chemfig{X_1 + 2 X_2}$)};
		%%%%%%
		\draw [{->[harpoon]}, transform canvas={ yshift=1.5pt}] (3.1,0.6)--(4.9,0.6) node [midway, above] {$\kk$};
		\draw [{[harpoon]<-}, transform canvas={ yshift=-1.5pt}] (3.1,0.6)--(4.9,0.6) node [midway, below] {$\kk$};
	}
where each box is a vertex of the graph; the top entry in each box is the stoichiometric complex of that vertex ($0$ and $\chemfig{X_1 + X_2}$), and the bottom entry in the parentheses is the kinetic-order complex ($\chemfig{2X_1}$ and $\chemfig{X_1 + 2 X_2}$). 
	The associated dynamical system of this generalized mass-action system is given by
	\eq{
		\frac{dx_1}{dt} &= \kk x_1^2 - \kk x_1 x_2^2, \\
		\frac{dx_2}{dt} &= \kk x_1^2 - \kk x_1 x_2^2.
	}
One can check that the set of vertex-balanced steady states is $Z_{\boldsymbol\kk} = \{ (t^2, t) : t >0 \}$. If $\bb x_0 = (0, \epsilon)^T$ where $0 < \epsilon < \frac{1}{4}$, then there are two vertex-balanced steady states in $\bb x_0 + S = \{ (r, \epsilon + r) : r \in \rr\}$. In particular, this implies that these vertex-balanced steady states cannot be globally stable in their stoichiometric compatibility class. \\

	Moreover, it is also possible for a unique vertex-balanced steady state (within its stoichiometric compatibility class) to be unstable. Consider the generalized mass-action system:
	\tikzc{
		\draw [rounded corners=5pt] (0,0) rectangle++(3,1.2);
		\draw [rounded corners=5pt] (5,0) rectangle++(3,1.2);
		\node at (1.5,0.85) {$
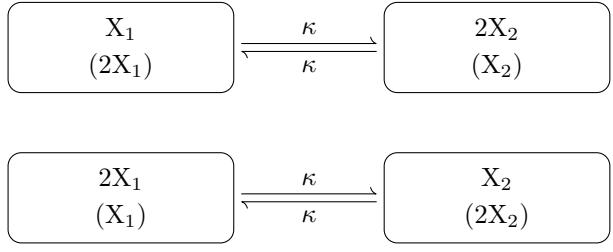
\chemfig{X_1}$};
		\node at (1.5,0.35) {($\chemfig{2X_1}$)};
		\node at (6.5,0.85) {$\chemfig{2 X_2}$};
		\node at (6.5,0.35) {($\chemfig{X_2}$)};
		%%%%%%
		\draw [{->[harpoon]}, transform canvas={ yshift=1.5pt}] (3.1,0.6)--(4.9,0.6) node [midway, above] {$\kk$};
		\draw [{[harpoon]<-}, transform canvas={ yshift=-1.5pt}] (3.1,0.6)--(4.9,0.6) node [midway, below] {$\kk$};
		%%%%%%%%%%%%%%%%%%%%%%%%%%%%%%%%%%%%%%%%%%%%%%%%%%%%%%%%%%%%%%
		\draw [rounded corners=5pt] (0,-2) rectangle++(3,1.2);
		\draw [rounded corners=5pt] (5,-2) rectangle++(3,1.2);
		\node at (1.5,-1.15) {$\chemfig{2 X_1}$};
		\node at (1.5,-1.65) {($\chemfig{X_1}$)};
		\node at (6.5, -1.15) {$\chemfig{X_2}$};
		\node at (6.5,-1.65) {($\chemfig{2 X_2}$)};
		%%%%%%
		\draw [{->[harpoon]}, transform canvas={ yshift=1.5pt}] (3.1,-1.4)--(4.9,-1.4) node [midway, above] {$\kk$};
		\draw [{[harpoon]<-}, transform canvas={ yshift=-1.5pt}] (3.1,-1.4)--(4.9,-1.4) node [midway, below] {$\kk$};
	}
		Its associated dynamical system is 
		\eq{
			\frac{dx_1}{dt} &= -\kk x_1^2 + \hphantom{2}\kk x_2 - 2 \kk x_1 + 2\kk x_2^2 , \\
			\frac{dx_2}{dt} &=\, \,2\kk x_1^2 - 2 \kk x_2 + \hphantom{2}\kk x_1 - \hphantom{2} \kk x_2^2.
		}	
This is an example of a reversible generalized mass-action system with $\delta_G = \tilde \delta_G = 0$, and its stoichiometric subspace $S$ and its kinetic-order subspace $\tilde S$ are $\rr^2$. There is a unique positive steady state $\boldsymbol{x}^* = (1,1)^T$, which is vertex-balanced; nonetheless, it can be shown that this steady state is a saddle point. Moreover, all solutions that start outside its stable manifold converge to the origin or infinity; in particular, the system is neither persistent nor permanent.

%%%%%%%%%%%%%%%%%%%%%%%%%%%%%%%%%%%%%%%%%%%%%%%%%%%%%%%%%%%%%%%%%%%%%%%%%%%%%%%%%%%%%%%%%%%%%%%%%%%%%%%%%%%%%%%%%%%%%%%
\section{An illustrative example}
\label{sec:MainEx}

We conclude by applying Theorem~\ref{thm:mainCRN} to the following example of a family of generalized mass-action systems. Let $a$, $b$, $\kk_i > 0$. Consider the generalized mass-action system $(G_{\boldsymbol\kk}, \Phi, \tilde \Phi)$
%%%%%%%%%%%%%%%%%%%%%%%%%%%%%%
	\tikzc{
		\draw [rounded corners=5pt] (0,0) rectangle++(3,1.2);
		\draw [rounded corners=5pt] (5,0) rectangle++(3,1.2);
		\draw [rounded corners=5pt] (2.5,-3) rectangle++(3,1.2);
		\node at (1.5,0.85) {$0$};
		\node at (1.5,0.35) {($0$)};
		\node at (6.5,0.85) {$
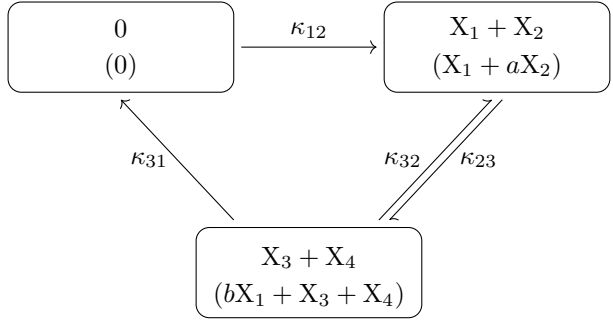
\chemfig{X_1 + X_2}$};
		\node at (6.5,0.35) {($\chemfig{X_1 +}a\chemfig{X_2}$)};
		\node at (4,-2.2) {$\chemfig{X_3 + X_4}$};
		\node at (4,-2.7) {($b\chemfig{X_1 + X_3}+ \chemfig{ X_4}$)};
		%%%%%%
		\draw [->] (3.1,0.6)--(4.9,0.6) node [midway, above] {$\kk_{12}$};
		\draw [->] (3,-1.7)--(1.5, -0.1) node [midway, left] {$\kk_{31}$};
		\draw [{->[harpoon]}, transform canvas={xshift=-2pt, yshift=0pt}] (5,-1.7)--(6.5, -0.1) node [midway, left] {$\kk_{32}$\,};
		\draw [{[harpoon]<-}, transform canvas={xshift=2pt, yshift=0pt}] (5,-1.7)--(6.5, -0.1) node [midway, right] {\,$\kk_{23}$};
	}
At each vertex (box), a stoichiometric complex (top entry) and a kinetic-order complex (second entry in parentheses) are assigned. Let $x_i$ be the concentration of species $\chemfig{X_i}$, for $1 \leq i \leq 4$, and $\bb x = (x_1,x_2,x_3,x_4)^T$. The stoichiometric complexes and kinetic-order complexes are
	\eq{
		\begin{array}{lclcl}
			\bb y_1 = (0,0,0,0)^T,
				&\quad &
				\bb y_2 = (1,1,0,0)^T,
				&\quad &
				\bb y_3 = (0,0,1,1)^T, \\
			\tilde{\bb y}_1 = (0,0,0,0)^T,
				&\quad &
				\tilde{\bb y}_2 = (1,a,0,0)^T,
				&\quad &
				\tilde{\bb y}_3 = (b,0,1,1)^T.
		\end{array}
	}
%$\bb y_1 = (0,0,0,0)^T$ and $\tilde{\bb y}_1 = (0,0,0,0)^T$; $\bb y_2 = (1,1,0,0)^T$ and $\tilde{\bb y}_2 = (1,a,0,0)^T$; and $\bb y_3 = (0,0,1,1)^T$ and $\tilde{\bb y}_3 = (b,0,1,1)^T$.
The associated dynamical system is
	\eq{
		\frac{d\bb x}{dt} &=
			\kk_{12} \bb x^{\tilde{\bb y}_1} \left(
				\bb y_2 - \bb y_1
			\right)
			+ \kk_{23} \bb x^{\tilde{\bb y}_2} \left(
				\bb y_3 - \bb y_2
			\right)
			+ \kk_{32} \bb x^{\tilde{\bb y}_3} \left(
				\bb y_2 - \bb y_3
			\right)
			+ \kk_{31} \bb x^{\tilde{\bb y}_3} \left(
				\bb y_1 - \bb y_3
			\right)
		\\&= \kk_{12}  %\left(
				\begin{pmatrix}
					1 \\ 1\\ 0 \\ 0
				\end{pmatrix} %y2
			%\right)
			+ \kk_{23} x_1x_2^a 
				\begin{pmatrix}
					-1 \\ -1\\ 1 \\ 1
				\end{pmatrix}
			+ \kk_{32} x_1^bx_3x_4 
				\begin{pmatrix}
					1 \\ 1\\ -1 \\ -1
				\end{pmatrix} 
			+ \kk_{31} x_1^bx_3x_4 %\left(
				 \begin{pmatrix}
					0\\0\\-1\\-1
				\end{pmatrix} %y3
			%\right)
			.
	}
Another way to write the system of differential equations is
	\eq{
		\frac{dx_1}{dt} &= \kk_{12} - \kk_{23} x_1x_2^a + \kk_{32} x_1^bx_3x_4, \\
		\frac{dx_2}{dt} &= \kk_{12} - \kk_{23} x_1x_2^a + \kk_{32} x_1^b x_3x_4, \\
		\frac{dx_3}{dt} &= \hspace{6ex} \kk_{23} x_1x_2^a - (\kk_{32} + \kk_{31} )x_1^bx_3x_4, \\
		\frac{dx_4}{dt} &= \hspace{6ex} \kk_{23} x_1x_2^a - (\kk_{32} + \kk_{31}) x_1^bx_3x_4.
	}
The stoichiometric subspace and the kinetic-order subspace are
	\eq{
		S = \Span_\rr \left\{
			\begin{pmatrix}
				1 \\ 1 \\ 0 \\ 0
			\end{pmatrix}
			,
			\begin{pmatrix}
				0 \\ 0 \\ 1 \\ 1
			\end{pmatrix}
		\right\},
		\qquad \text{and} \qquad
		\tilde S = \Span_\rr\left\{
			\begin{pmatrix}
				1 \\ a \\ 0 \\ 0
			\end{pmatrix}
			,
			\begin{pmatrix}
				b \\ 0 \\ 1 \\ 1
			\end{pmatrix}
		\right\}
	}
respectively. Their sign vectors are
	\eq{
		\sigma(S) &=
			\left\{
				\begin{pmatrix}
				 0 \\ 0 \\ 0 \\ 0
				\end{pmatrix},
				\begin{pmatrix}
				 + \\ + \\ + \\ +
				\end{pmatrix},
				\begin{pmatrix}
				 0 \\ 0 \\ + \\ +
				\end{pmatrix},
				\begin{pmatrix}
				 - \\ - \\ + \\ +
				\end{pmatrix},
				\begin{pmatrix}
				 - \\ - \\ 0 \\ 0
				\end{pmatrix},
				\,\cdots\,
			\right\},
\intertext{and}
		\sigma(\tilde S) &=
			\left\{
				\begin{pmatrix}
				 0 \\ 0 \\ 0 \\ 0
				\end{pmatrix},
				\begin{pmatrix}
				 + \\ + \\ + \\ +
				\end{pmatrix},
				\begin{pmatrix}
				 + \\ 0 \\ + \\ +
				\end{pmatrix},
				\begin{pmatrix}
				 + \\ - \\ + \\ +
				\end{pmatrix},
				\begin{pmatrix}
				 0 \\ - \\ + \\ +
				\end{pmatrix},
				\begin{pmatrix}
				 - \\ - \\ + \\ +
				\end{pmatrix},
				\begin{pmatrix}
				 - \\ - \\ 0 \\ 0
				\end{pmatrix},
				\, \cdots\,
			\right\},
	}
where the dots indicate the negatives of the listed sign vectors. By visual inspection, we find that $\sigma(S) \subseteq \cl{\sigma(\tilde S)}$. Moreover, one can check that the deficiency $\delta_G$ and the kinetic-order deficiency $\tilde\delta_G$ are zero. Therefore, Corollary~\ref{cor:last} applies and we conclude that, for any choice of rate constants, every stoichiometric compatibility class contains exactly one positive steady state, which is vertex-balanced.

\section*{Acknowledgments}

This project started during the AIM SQuaRE workshop on ``Dynamical properties of deterministic and stochastic models of reaction networks''. The authors thank the American Institute of Mathematics for the productive environment and hospitality. GC was partially supported by NSF-DMS grants 1412643 and 1816238; SM was supported by the Austria Science Fund (FWF), project P28406; CP was partially supported by NSF-DMS grant 1517577; and PYY was partially supported by a NSERC PGS-D award and NSF-DMS grant 1412643.

\bibliographystyle{amsxport}
\bibliography{cit}

%%% If direct bibliography entries. Lists all in order of entry
%\begin{bibdiv}
%\begin{biblist}
%	\bib{Sokal96}{article}{
%		title={Trangressing the boundaries},
%		subtitle={Toward a transformative hermeneutics of quantum gravity},
%		author={Sokal, Alan},
%		journal={Social Text},
%		volume={46/47},
%		date={1996},
%		pages={217--252}
%	}
%\end{biblist}
%\end{bibdiv}

\end{document}